\documentclass[11pt]{article}

\usepackage[dvips]{epsfig}
\usepackage{amssymb}

\DeclareMathAlphabet{\mathgo}{U}{euf}{m}{n}

\def\Fix{{\rm Fix}}
\newcommand{\qbinomial}[3]{\left[{\displaystyle #1} \atop {\displaystyle #2}\right]_#3}

\newtheorem{theo}{Theorem}
\newtheorem{prop}[theo]{Proposition} 
\newtheorem{defin}{Definition}
\newtheorem{examp}{Example} 
\newtheorem{cor}[theo]{Corollary}

\newtheorem{observat}[theo]{Remark}
\newtheorem{lemma}[theo]{Lemma}

\newenvironment{rem*}{
\begin{trivlist}\item
\textbf{Remark}\hspace{\blspace}\rm}
{\end{trivlist}}

\newlength{\maxitemlgth} 
\newcounter{desc}
\newlength{\blspace} 
\newsavebox{\qedbox} 

\sbox{\qedbox}{%
\begin{picture}(5.1,5)%
\put(0,0){\framebox(5,5){}}%
\end{picture}} %
\newcommand{\qed}{\hfill\usebox{\qedbox}}%
\newenvironment{proof}[1][]%
{\begin{trivlist}\item{\textbf{Proof#1.}\hspace{\blspace}}}%
{\qed\end{trivlist}}

\begin{document}

\title{Enumeration of Symmetry Classes of Convex Polyominoes in the Square Lattice\thanks{Work partially supported by NSERC (Canada) and FCAR (Qu\'ebec).}}
\author{\textsc{P. Leroux, E. Rassart, A. Robitaille}}
\maketitle
\centerline{\textsl{LaCIM, Department of Mathematics}}
\centerline{\textsl{Universit\'e du Qu\'ebec \`a Montr\'eal} (UQAM)}

\begin{abstract}
{This paper concerns the enumeration of rotation-type and congruence-type convex polyominoes on the square lattice. These can be defined as orbits of the groups $\mathfrak{C}_4$, of rotations, and $\mathfrak{D}_4$, of symmetries, of the square, acting on (translation-type) polyominoes. In virtue of Burnside's Lemma, it is sufficient to enumerate the various symmetry classes (fixed points) of polyominoes defined by the elements of $\mathfrak{C}_4$ and $\mathfrak{D}_4$. Using the Temperley--Bousquet-M\'elou methodology, we solve this problem and provide explicit or recursive formulas for their generating functions according to width, height and area. We also enumerate the class of asymmetric convex polyominoes, using M\"obius inversion, and prove that their number is asymptotically equivalent to the number of convex polyominoes, a fact which is empirically evident.}
\end{abstract}

\renewcommand\abstractname{\protect R\'esum\'e}

\begin{abstract}
{Cet article porte sur l'\'enum\'eration de polyominos \`a rotations et r\'eflexions pr\`es dans un r\'eseau carr\'e. Ces polyominos peuvent \^etre consid\'er\'es comme des orbites dans l'action du groupe di\'edral $\mathfrak{D}_4$ du carr\'e et de son sous-groupe $\mathfrak{C}_4$, le groupe des rotations du carr\'e, sur l'ensemble de polyominos convexes (\`a translations pr\`es). Le lemme de Burnside r\'eduit alors le probl\`eme \`a l'\'enum\'eration des ensembles de points fixes (classes de sym\'etrie) de l'action. Utilisant la m\'ethodologie de Temperley--Bousquet-M\'elou, nous donnons des formules explicites ou r\'ecursives pour les s\'eries g\'en\'eratrices de ces types de polyominos selon la largeur, la hauteur et l'aire. Nous \'enum\'erons aussi, \`a l'aide de l'inversion de M\"obius, la classe des polyominos convexes asym\'etriques et prouvons que ce nombre est asymptotique \`a celui des polyominos convexes, un fait mis en \'evidence exp\'erimentalement.}
\end{abstract}

%
%
\section{Introduction}

The enumeration of polyominoes and of self-avoiding walks of a square lattice is a deep and difficult combinatorial problem. These structures are important in statistical physics, as discrete models for polymers, clusters, percolation and phase transitions. See, for instance, Privman and \u{S}vraki\'c \cite{VPNMS}, Viennot \cite{XGV}, Bousquet-M\'elou \cite{MBMb}, Guttmann, Prellberg and Owczarek \cite{AJGTPALO}, Langlands, Pouliot and Saint-Aubin \cite{RLPPYSA}, Prellberg and Brak \cite{TPRB} and Hughes \cite {BDH1} and \cite{BDH2}.

\begin{figure}[!ht]
\begin{center}
\includegraphics{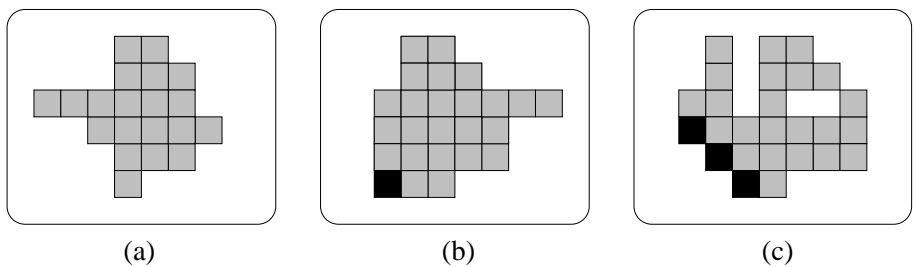}
\caption{\small Polyominoes: (a) convex; (b) directed (one source) convex; (c) directed (diagonal source).}
\label{PolyoDef}
\end{center}
\end{figure}

Recent research in this area has mainly dealt with polyominoes up to translation (where two polyominoes are equivalent if some translation of the plane maps one onto the other). The equivalence classes are then called \emph{translation-type} polyominoes, or simply polyominoes when the context is clear. The exact enumeration of polyominoes is an open problem, but considerable progress has been made for classes of polyominoes defined by special properties of convexity or directedness (see \cite{MBMc}). A polyomino is said to be \textsl{convex} when its intersection with any vertical or horizontal line is connected. A polyomino $P$ is (North-East) \textsl{directed} if each of its cells can be reached from some distinguished cells (the source cells) by a path that is contained in $P$ and has only North and East steps (see Figure~\ref{PolyoDef}). Formulas have been published for the generating function of convex polyominoes according to perimeter (\cite{MPDGV}, \cite{DK} and \cite{AJGIGE}), to height and width (\cite{KYLSJC}) and joint parameters area, height and width (\cite{KYL}, \cite{MBMJMF} and \cite{MBMa}).

In this paper, we consider \emph{rotation-type} polyominoes, that is equivalence classes of polyominoes under rotations, and \emph{congruence-type} polyominoes, where rotations and reflections are permitted. They occur naturally as pieces that can freely move in space, as used in plane packing problems (see \cite{DAK} and \cite{SG}). For example, there are 19 (translation-type) tetrominoes, 7 rotation-type tetrominoes and 5 congruence-type tetrominoes (tetrominoes are polyominoes with area 4).

\begin{figure}[!ht]
\begin{center}
\includegraphics{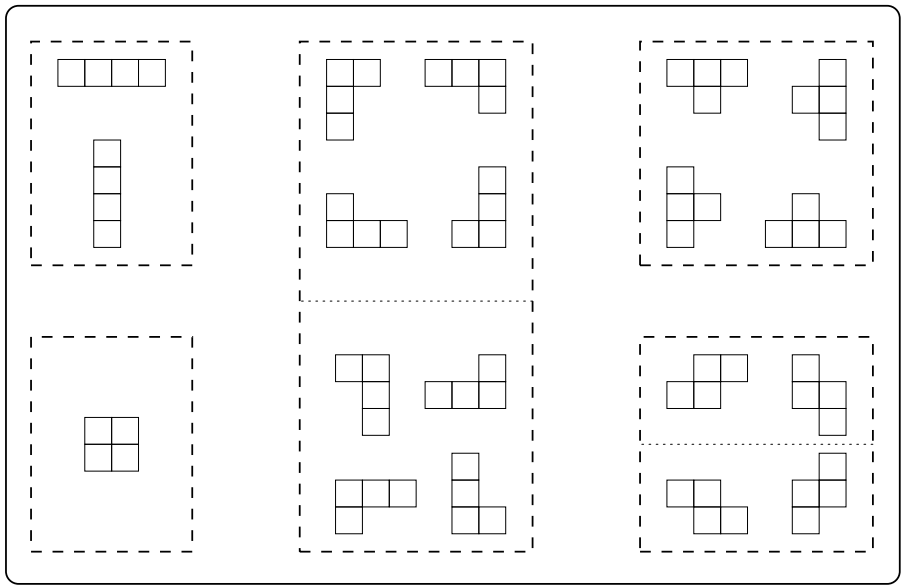}
\caption{\small The 19 tetrominoes and their equivalence classes.}
\label{Tetraminoes}
\end{center}
\end{figure}

\begin{figure}[!ht]
\begin{center}
\includegraphics{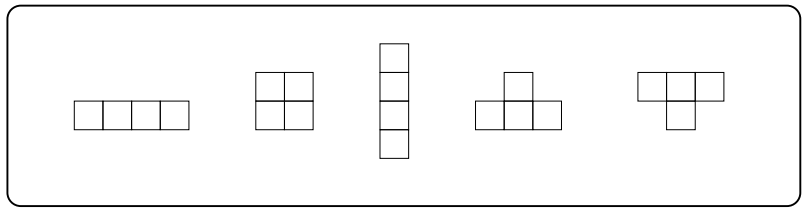}
\caption{\small Vertically symmetric tetrominoes.}
\label{VertSymTetra}
\end{center}
\end{figure}

As Figure~\ref{Tetraminoes} illustrates, congruence-type polyominoes can be seen as orbits of the action of the group $\mathfrak{D}_4$ of symmetries of the square on the set of translation-type polyominoes. Similarly, rotation-type polyominoes are orbits of the subgroup $\mathfrak{C}_4$ of $\mathfrak{D}_4$ containing the rotations. In virtue of Burnside's Lemma, the enumeration of rotation-type or congruence-type polyominoes reduces to the enumeration of the various symmetry classes, i.e. the sets of fixed polyominoes determined by the group elements of $\mathfrak{C}_4$ and $\mathfrak{D}_4$. Figure~\ref{VertSymTetra} illustrates the class of vertically symmetric tetrominoes, i.e. tetrominoes left fixed by a reflection about a vertical axis.

The problem of enumerating rotation-type and congruence-type polyominoes is mentionned in Golomb \cite{SG} and Klarner \cite{DAK}, but only experimental results, for polyominoes of area up to $24$, are reported, together with asymptotic estimates. Here we solve this problem for the class of convex polyominoes. The study of convex polyominoes was suggested by D.~Knuth (see \cite{DAKRLR}). Convex polyominoes also appear naturally in statistical mechanics as self-avoiding polygons whose number of steps is equal to the perimeter of their minimal bounding rectangle (see \cite{AJGIGE}). Using elementary methods and the Temperley--Bousquet-M\'elou methodology \cite{MBMa}, we give explicit or recursive formulas for the generating functions according to perimeter and area of each of their symmetry classes. These results are then combined to enumerate rotation-type and congruence-type convex polyominoes. Using M\"obius inversion in the lattice of subgroups of $\mathfrak{D}_4$, we also count asymmetric convex polyominoes. Since the symmetry class defined by each subgroup of $\mathfrak{D}_4$ (except the trivial one) is significantly smaller than the whole class, we conclude that the number of asymmetric convex polyominoes, as the area goes to infinity, is asymptotic to the number of convex polyominoes. In other words, almost all convex polyominoes are asymmetric.

The paper is organized as follows. In section 2, we review some basic results that are used in the sequel on group actions and convex polyomino enumeration, according to width, height and area. Polyomino classes that are covered include partitions, stacks, in particular the new class of shifted stacks, directed convex and convex polyominoes. We then treat successively the symmetry classes of convex polyominoes in sections 3 (rotational symmetries), 4 (horizontal and vertical symmetries) and 5 (diagonal symmetries). The enumeration of diagonally symmetric convex polyominoes requires the determination of a subsidiary class of diagonal source (shifted) directed convex polyominoes, in subsection 5.2. These results are then combined to enumerate rotation-type polyominoes (subsection 3.3) and congruence-type polyominoes (section 6) according to half-perimeter and area. Finally section 7 is devoted to asymmetric convex polyominoes. Some numerical results are found in Appendices A and B.

\vspace{5mm}

\centerline{\textbf{Acknowledgements}}
\vspace{5mm}
The first author would like to thank Mireille Bousquet-M\'elou, Richard Brak, Dominique Gouyou-Beauchamps and Tony Guttmann for stimulating and useful discussions.

%
%
\section{Preliminaries}

%
%
\subsection{Review of convex polyomino enumeration}
\label{sec:review}

We recall here some basic and known results on the enumeration of some classes of translation-type polyominoes that will be used in this paper, namely Ferrers diagrams or partition polyominoes ($P$), partitions with empty parts allowed ($P_0$), (shifted) partitions with distinct parts ($P_S$), varieties of stack or pyramid polyominoes ($T$, $T_0$ and $T_S$), directed (one source) convex ($D$), and convex ($C$) polyominoes. For each of these classes, we give the generating series $F(x,y,q)$ where the variable $x$ marks the width of the polyominoes, $y$ marks the height and $q$ the area. For example, the weight of the convex polyomino of Figure~\ref{PolyoDef}(a) according to these parameters is $x^7y^6q^{20}$. We will often use the same letters $P$, $T$, $D$, \ldots to denote a class of polyominoes, its generating function, or a generic polyomino in the class.

\subsubsection{Partition polyominoes}
\label{sec:partition}

\begin{figure}[ht]
\begin{center}
\includegraphics{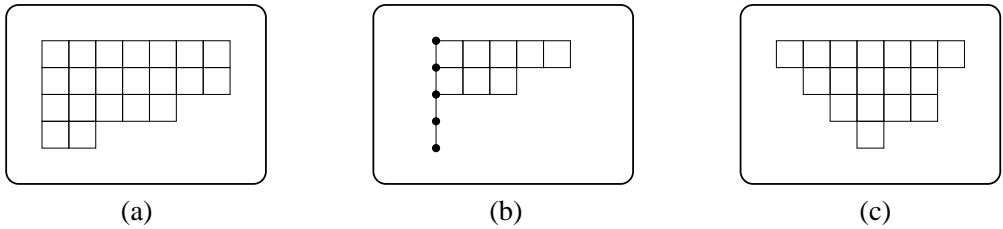}
\caption{\small Partition polyominoes. (a) $P$: general; (b) $P_0$: with empty parts allowed; (c) $P_S$: shifted.}
\label{PartDef}
\end{center}
\end{figure}

A partition polyomino is the Ferrers diagram of a non-empty partition (see Figure~\ref{PartDef}(a)). It can be viewed as a left justified pile of weakly decreasing rows (from top to bottom). Let $P(x,y,q)$ denote their generating series. In this context, the variable $x$ marks the size of the maximum part and $y$ the number of parts. Let also $P_0(x,y,q)$ denote the generating series of partitions where empty rows (contributing to the height) are permitted (see Figure~\ref{PartDef}(b)), and $P_S(x,y,q)$, of partitions with distinct parts (diagonally shifted) (see Figure~\ref{PartDef}(c)).

We then have the following classical results (see \cite{GEA}), for which direct counting arguments can be given:

\begin{eqnarray}
P(x,y,q) & = & \sum_{m \geq 1}\frac{x^myq^m}{(yq)_m} \nonumber \\
         & = & \sum_{m, k\, \geq 1}x^my^kq^{m+k-1}\qbinomial{m+k-2}{k-1}{q}\,,\\
         &   & \nonumber \\
P_0(x,y,q) & = & \sum_{m \geq 0}\frac{x^myq^m}{(y)_{m+1}} \nonumber \\
           & = & \sum_{m \geq 0, k \geq 1}x^my^kq^m\qbinomial{m+k-1}{k-1}{q}\,,\\
           &   & \nonumber \\
P_S(x,y,q) & = & \sum_{m \geq 1}x^myq^m(-yq)_{m-1} \nonumber \\
           & = & \sum_{m \geq 1}\sum_{k = 1}^mx^my^kq^{m + \left(k \atop 2\right)}\qbinomial{m-1}{k-1}{q}\,,
\end{eqnarray}
where we have used the familiar notation
\begin{equation}
(a)_n = (a; q)_n = (1-a)(1-aq)\ldots (1-aq^{n-1})\,,  \textrm{for $n \geq 0$}
\end{equation}
so that in particular $(-yq)_m = (1+yq)(1+yq^2)\ldots (1+yq^m)$, and where, for $0 \leq k \leq n$,
\begin{equation}
\qbinomial{n}{k}{q} = \frac{(q)_n}{(q)_k(q)_{n-k}}
\end{equation}
is the $q$-binomial (Gaussian) polynomial.

\subsubsection{Stack polyominoes}
\label{sec:stack}

\begin{figure}[ht]
\begin{center}
\includegraphics{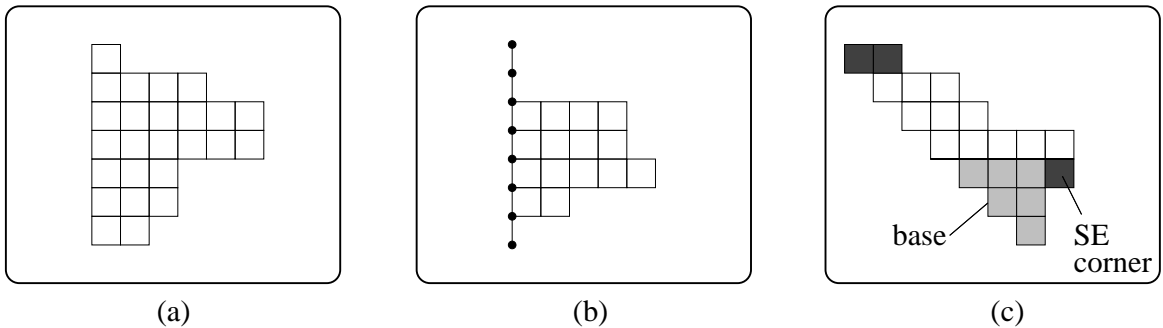}
\caption{\small Stack polyominoes. (a) $T$: general; (b) $T_0$: with empty parts allowed; (c) $T_S$: shifted stacks.}
\label{StackDef}
\end{center}
\end{figure}

A \textsl{stack} or \textsl{pyramid} polyomino is defined as a left justified unimodal (i.e. weakly increasing, then weakly decreasing) pile of rows (see Figure~\ref{StackDef}(a)). Let $T(x,y,q)$ ($T$ for French ``tas'') denote the generating series of stacks, and $T_0(x,y,q)$, the generating series of stacks with empty rows allowed at either end (see Figure~\ref{StackDef}(b)). 

By extracting coefficients, define, for $n \geq 1$,
\begin{equation}
T_n(x,q) = [y^n]\,T(x,y,q) \quad \textrm{and} \quad T_{0,n}(x,q) = [y^n]\,T_0(x,y,q)\,,
\end{equation}
as the width and area generating series for stack polyominoes having height equal to $n$. It is  clear that $T(x,y,q) = xT_0(x,yq,q)$ and that $T_n(x,q) = xq^nT_{0,n}(x,q)$. Hence any enumerative result for the class $T_0$ will yield a result for $T$. By concentrating on the width of the largest row, we obtain
\begin{equation}
\label{eqn:stack1}
T_0(x,y,q) = \sum_{m \geq 0}\frac{x^myq^m}{(y)_{m+1}(y)_{m}}
\end{equation}
and 
\begin{equation}
\label{eqn:stacks}
T_{0,n}(x,q) = \sum_{m \geq 1}\sum_{j=1}^{n}x^mq^m\qbinomial{m+j-1}{m}{q}\qbinomial{m+n-j-1}{m-1}{q}\,.
\end{equation}

The following simple expression is due to Bousquet-M\'elou and Viennot \cite{MBMXGV}.

\begin{prop}The width and area generating series for stacks having $n$ (possibly empty) rows is given by
\begin{equation}
\label{eqn:stack2}
T_{0,n}(x,q) = \frac{V_n}{(xq)_n}\,,\quad n \geq 1\,,
\end{equation}
where $V_n = V_n(x,q)$ are polynomials in $x$ and $q$ defined by the recurrence
\begin{equation}
\label{Vnrec}
V_n = 2V_{n-1} + (xq^{n-1} - 1)V_{n-2}\ , \quad n \geq 2\,,
\end{equation}
with $V_0 = V_1 = 1$. Moreover, $V_n$ can be given explicitly as 
\begin{equation}
\label{Vnclosed}
V_n = 1 + \sum_{1 \leq k \leq n/2}\left(x^kq^{k^2}\sum_{m=k}^{n-k}\qbinomial{m}{k}{q}\qbinomial{n-m-1}{k-1}{q}\right)\,.
\end{equation}
\end{prop}

\begin{proof} Setting $T_{0,0}(x,q) = 1$ by convention (the empty stack), it is easy to see that for $n \geq 2$,
\begin{equation}
T_{0,n} = xq^nT_{0,n} + 2T_{0,n-1} - T_{0,n-2}\,,
\end{equation}
so that $(1-xq^n)T_{0,n} = 2T_{0,n-1} - T_{0,n-2}$. Since $T_{0,1} = 1/(1-xq) = 1/(xq)_1$, we see that indeed $T_{0,n} = V_n/(xq)_n$ where $V_n$ satisfies $V_0 = V_1 = 1$ and recurrence~(\ref{Vnrec}). This recurrence can also be used to prove that for $n \geq 1$, $V_n$ is the generating polynomial for $T_0$-stacks of height $n-1$ having the additional property that any column exceeds strictly its right neighbor both at the top and at the bottom. Formula~(\ref{Vnclosed}) then follows easily.
\end{proof}

\subsubsection{Shifted stacks}
\label{sec:diagrowstack}

In analogy with shifted partitions $P_S$, we introduce the class $T_S$ of (diagonally) shifted stacks (Figure~\ref{StackDef}(c)), and we denote by $T_S(u,x,y,q)$ their generating series, where $u$ is an extra variable which marks the top row width. For example, the shifted stack of Figure~\ref{StackDef}(c) has global weight $u^2x^8y^7q^{20}$. This generating series will be used in Section~\ref{sec:diagsymclass}.

Define the \textsl{South-East corner} of a shifted stack $T_S$ to be the bottommost cell of the rightmost column of $T_S$. The \textsl{base} of $T_S$ consists of all rows below (and including) the South-East corner row. It is a shifted partition $P_S$ with the additional property that the first row exceeds the other rows, i.e. that the first part of the partition exceeds the others by at least $2$. Let $P_1(u,x,y,q)$ denote the generating series of all such possible bases of shifted stacks. We obviously have
\begin{equation}
P_1(u,x,y,q) = uxyq + y\sum_{k \geq 2}(uxq)^k(-yq)_{k-2}\,.
\end{equation}

\begin{prop}The generating function $T_S(u) = T_S(u,x,y,q)$ of shifted stacks satisfies the functional equation
\begin{equation}
\label{Rfunceqn}
T_S(u) = P_1(u,x,y,q) + \frac{xyu^2q^2}{1-uq}\left(T_S(1) - T_S(uq)\right)
\end{equation}
\end{prop}

\begin{proof}Equation~(\ref{Rfunceqn}) states that a shifted stack is either a $P_1$-polyomino or is obtained from a $P_1$-polyomino by iterating the procedure that adds a row on top of a shifted stack. At each step, the added row should be shifted one cell to the left and should not exceed the previous row on the right hand side. Now the term $xyu^2q^2(1-uq)^{-1}T_S(1)$ enumerates the cases where a shifted row of any length ($\geq 2$) is added on top of a shifted stack, and $xyu^2q^2(1-uq)^{-1}T_S(uq)$ counts the cases where the added row does exceed on the right, so that this term should in fact be subtracted. 
\end{proof}

\begin{figure}[!ht]
\begin{center}
\includegraphics{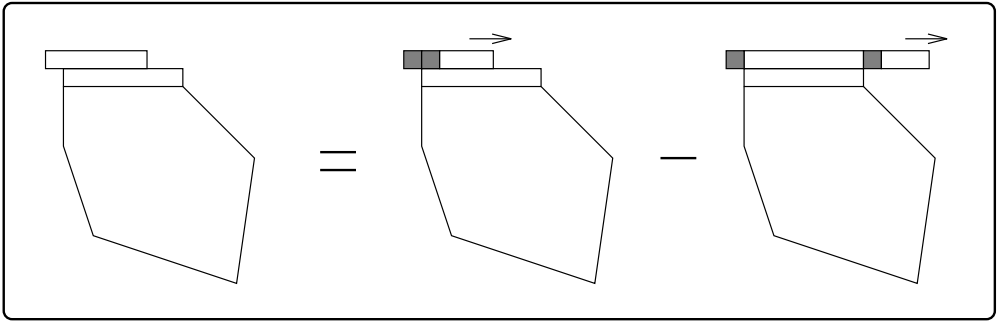}
\caption{A variation of construction B of Bousquet-M\'elou.}
\label{DiagRec}
\end{center}
\end{figure}

The above procedure of adding a left shifted non-exceeding row (see Figure~\ref{DiagRec}) can be seen as a variation of construction B of Bousquet-M\'elou \cite{MBMa}, where Temperley's methodology is translated into linear $q$-equations. (see also \cite{SFDS}). Using Lemma 2.3 of \cite{MBMa} which solves~(\ref{Rfunceqn}) by iteration, we find the following closed form:

\begin{cor}The generating series $T_S(u) = T_S(u,x,y,q)$ for shifted stacks is given by
\begin{equation}
\label{funceqnsol}
T_S(u) = \frac{E(u) + E(1)F(u) - E(u)F(1)}{1 - F(1)}\,,
\end{equation}
where
\begin{eqnarray}
E(u) & = & \sum_{n \geq 0}\frac{(-1)^nx^{n+1}y^nu^{2n}q^{n^2+n}P_1(uq^n,x,y,q)}{(uq)_n}\,, \\
F(u) & = & \sum_{n \geq 0}\frac{(-1)^nx^{n+1}y^{n+1}u^{2n+2}q^{n^2+3n+2}}{(uq)_{n+1}}\,.
\end{eqnarray}\hfill $\square$
\end{cor}

Note that the functional equation~(\ref{Rfunceqn}) and its solution~(\ref{funceqnsol}) remain valid for other possible classes $P_1$ of polyominoes as bases. The method can then be adapted to enumerate Dyck words ($P_1$ is a single cell) and left factors of Dyck words ($P_1$ is a single row) according to length, area and size of the first stretch of $x$'s (marked by $u$).

\subsubsection{Directed convex and convex polyominoes}
\label{sec:dirconv}

We now give without proof the generating function $D(s,x,y,q)$ for directed convex polyominoes, where the variable $s$ marks the leftmost column height, borrowed from Bousquet-M\'elou \cite{MBMa} (Theorem 3.4, with $t = 1$).

\begin{prop}The generating function for directed convex polyominoes is given by
\begin{equation}
D(s,x,y,q) = y\frac{M_1(s)J_0(1) - M_1(1)J_0(s) + M_1(1)}{J_0(1)}\,,
\end{equation}
where
\begin{equation}
\label{J0eqn}
J_0(s) = \sum_{n \geq 0}\frac{(-1)^nx^ns^nq^{\left(n+1 \atop 2\right)}}{(sq)_n(syq)_n}
\end{equation}
and
\begin{equation}
\label{M1eqn}
M_1(s) = \sum_{n \geq 1}\left(\frac{x^nq^n}{(syq)_n}\left(\frac{s}{(syq)_{n-1}} + \sum_{m=1}^{n-1}\frac{(-1)^ms^mq^{\left(m \atop 2\right)}}{(sq)_{n-1}(syq^{m+1})_{n-m-1}}\right)\right)\,.
\end{equation}\hfill $\square$
\end{prop}

We also extract from \cite{MBMa} (Theorem 4.4, with $t = 1$) the generating function $C(x,y,q)$ for convex polyominoes. We introduce the series $E(s) = E(s,x,y,q)$ by
\begin{equation}
\label{Eeqn}
E(s) = 2y^2\frac{M_1(1)}{J_0(1)}\alpha(s) - ya(s)\,,
\end{equation}
where $J_0$ and $M_1$ are given by~(\ref{J0eqn}) and~(\ref{M1eqn}),
\begin{equation}
\alpha(s) = \sum_{m \geq 1}\left(\frac{(-1)^mx^ms^mq^{\left(m+1 \atop 2\right)}}{(sq)_m}\sum_{n=1}^{m}\frac{(-1)^ns^nq^{\left(n+1 \atop 2\right)}}{(sq)_{n-1}(syq^n)_{m-n+1}}\right)
\end{equation}
and
\begin{eqnarray}
a(s) & = & \sum_{m \geq 1}\left(\frac{x^mq^m}{1-syq^m}\left(-\frac{s^{2m-1}q^{m(m-1)}}{(sq)^2_{m-1}} + y\sum_{n-1}^{m-1}\frac{s^{2n}q^{n^2}(syq^n-2)}{(sq)^2_{n-1}(syq^n)^2_{m-n}}\right.\right. \nonumber \\
 &  & \qquad + \left.\left.2y\sum_{1 \leq n \leq k < m}\frac{(-1)^{n+k}s^{n+k}q^{\left(n+1 \atop 2\right)}q^{\left(k \atop 2\right)}}{(sq)_{n-1}(syq^n)_{m-n}(sq)_k(syq^{k+1})_{m-k-1}}\right)\right)\,.
\end{eqnarray}

\begin{prop}The generating series $C(x,y,q)$ for convex polyominoes is given by
\begin{equation}
C(x,y,q) = J_1E'(1) - K_1E(1)\,,
\end{equation}
where $E(s)$ is given by~(\ref{Eeqn}), the derivative $E'$ is taken with respect to $s$,
\begin{equation}
J_1 = \sum_{n \geq 1}\frac{x^nq^n}{(q)_{n-1}(q)_n}\,,
\end{equation}
and
\begin{equation}
K_1 = -1 + \sum_{n \geq 1}\left(\frac{x^nq^n}{(q)_{n-1}(q)_n}\left(\sum_{k=1}^{n-1}\frac{2q^k}{1-q^k}+\frac{q^n}{1-q^n}\right)\right)\,.
\end{equation}\hfill $\square$
\end{prop}

%
%
\subsection{The action of the dihedral group $\mathfrak{D}_4$ on polyominoes}
\label{sec:mobius}

The dihedral group $\mathfrak{D}_4$ is the symmetry group of the square. It has eight elements, usually represented as $1$, $r$, $r^2$, $r^3$, $h$, $v$, $d_1$ and $d_2$, where $1$ denotes the identity element, $r$ denotes a rotation by a right angle, $h$ and $v$ reflections with respect to the horizontal and vertical axes respectively, and $d_1$ and $d_2$ reflections about the two diagonal axes of the square. The (nontrivial) subgroups of $\mathfrak{D}_4$ are the following:

\begin{displaymath}\left.
\begin{array}{cll}
\bullet & \textrm{The cyclic subgroups of order 2:} & \langle r^2\rangle, \langle h\rangle, \langle v\rangle, \langle d_1\rangle, \langle d_2\rangle \\
\bullet & \textrm{The cyclic subgroup of order 4:} & \mathfrak{C}_4\ =\ \langle r\rangle\ =\ \{1, r, r^2, r^3\} \\
\bullet & \textrm{The (non-cyclic) subgroups of order 4:} & \langle h, v\rangle\ =\ \{1, h, v, r^2\}\\
& & \langle d_1, d_2\rangle\ =\ \{1, d_1, d_2, r^2\}
\end{array}\right.
\end{displaymath} 

The group $\mathfrak{D}_4$ acts on (translation-type) polyominoes in a natural way, by rotating or reflecting them. Let $\mathcal{P}$ be a finite set of polyominoes closed under the action of $\mathfrak{D}_4$. We could take, for example, the set of convex polyominoes with a given area and/or perimeter. Then any subgroup $G$ of $\mathfrak{D}_4$ acts on $\mathcal{P}$ and the $G$-orbits are the equivalence classes of polyominoes modulo the transformations in $G$. The set of $G$-orbits of polyominoes in $\mathcal{P}$ is denoted by $\mathcal{P}/G$. In particular, the set of rotation-type polyominoes is precisely $\mathcal{P}/\mathfrak{C}_4$, and the set of congruence-type polyominoes, $\mathcal{P}/\mathfrak{D}_4$.

The number of distinct $G$-orbits of polyominoes in $\mathcal{P}$ is given by the Cauchy-Frobenius formula (alias Burnside's Lemma):
\begin{equation}
|\mathcal{P}/G| = \frac{1}{|G|}\sum_{g \in G}|\Fix_\mathcal{P}(g)|\,,
\label{BurnsideL}
\end{equation}
where $\Fix_\mathcal{P}(g)$ is the subset of $\mathcal{P}$ consisting of the invariant polyominoes under the transformation $g$. The $\Fix_\mathcal{P}(g)$ is called the \emph{symmetry class} of $g$, and we call its elements \textsl{$g$-symmetric polyominoes}. Formula~(\ref{BurnsideL}) extends to infinite classes of polyominoes by taking a weighted cardinality, e.g. the generating series with respect to half-perimeter (variable $t$) and area (variable $q$).

In the following sections we compute the generating series of the symmetry classes $\Fix_\mathcal{P}(g)$ for the set $C$ of all convex (translation-type) polyominoes, and for $g \in \mathfrak{D}_4$. We use the notation
$$F_g(t,q) = \textrm{generating series of}\ \Fix_{C}(g)\,.$$

Setting
$$\big(C/G\big)(t,q) = \textrm{generating series of $G$-orbits of convex polyominoes}$$
we have
\begin{equation}
\big(C/G\big)(t,q) = \frac{1}{|G|}\sum_{g \in G}F_g(t,q)\,.
\end{equation}

The concept of symmetry classes of polyominoes is extended to any subgroup $G$ of $\mathfrak{D}_4$ by setting
\begin{eqnarray*}
S_{\geq G} & = & \cap_{g \in G}\Fix_\mathcal{P}(g) = \{P \in \mathcal{P}\  |\  \textrm{Stab}(P) \supseteq G\}\,,\\
S_{= G} & = & \{P \in \mathcal{P}\  |\  \textrm{Stab}(P) = G\}\,,
\end{eqnarray*}
where Stab($P$) denotes the stabilizer of $P$, that is the subgroup of $\mathfrak{D}_4$ containing the transformations which leave $P$ unchanged. Observe in particular that $S_{= 0}$ is the set of asymmetric polyominoes. We let $F_{\geq G}$ and $F_{= G}$ denote the generating series of the sets $S_{\geq G}$ and $S_{= G}$ respectively. The series $F_{\geq G}$ for the subgroups $G$ of $\mathfrak{D}_4$ are easier to compute directly than $F_{= G}$. We obviously have
\begin{equation}
F_{\geq G} = \sum_{H \supseteq G}F_{= H}\,,
\end{equation}
where $H$ ranges over the lattice of subgroups of $\mathfrak{D}_4$, denoted by $\mathcal{S}(\mathfrak{D}_4)$.

Using M\"obius inversion (see Rota \cite{GCR}) in the lattice $\mathcal{S}(\mathfrak{D}_4)$, we then find
\begin{equation}
F_{= G} = \sum_{H \supseteq G}\mu(G, H)\,F_{\geq H}\,,
\end{equation}
where $\mu(G, H)$ is the value of the M\"obius function in $\mathcal{S}(\mathfrak{D}_4)$. In particular, taking $G = 0 = \{1\}$, the minimum element in $\mathcal{S}(\mathfrak{D}_4)$, we obtain the generating series for the asymmetric convex polyominoes:
\begin{equation}
\label{eqn:assym}
F_{= 0} = \sum_{H \in \mathcal{S}(\mathfrak{D}_4)}\mu(0, H)\,F_{\geq H}\,.
\end{equation}
\begin{figure}[!ht]
\begin{center}
\includegraphics{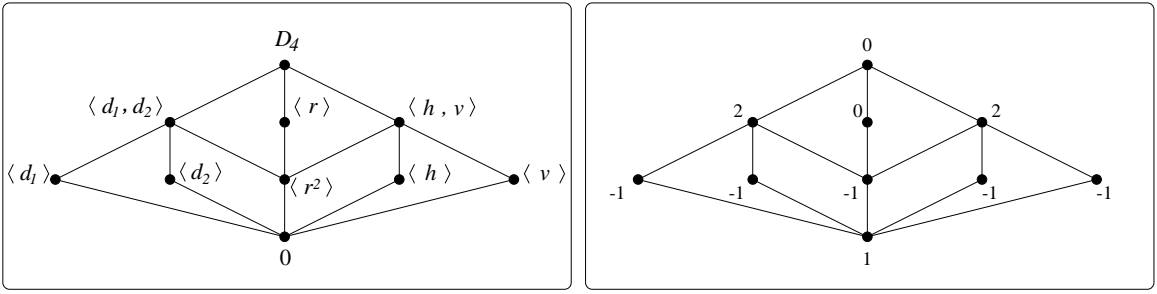}
\caption{\small Lattice of the subgroups of $\mathfrak{D}_4$ and the M\"obius function $\mu(0, H)$ on the lattice.}
\label{TreillisIntro}
\end{center}
\end{figure}

The value of the M\"obius function $\mu(0, H)$ is easy to compute recursively (see Rota \cite{GCR}), using the formulas
\begin{equation}
\mu(0, y) = - \sum_{x < y}\mu(0, x)\ ,\qquad \mu(0, 0) = 1\,.
\end{equation}
These values of $\mu$ on $\mathcal{S}(\mathfrak{D}_4)$ are given in Figure~\ref{TreillisIntro}.

Note that for a cyclic subgroup $\langle g\rangle$ of $\mathfrak{D}_4$, $F_{\geq \langle g\rangle}$ is precisely $F_g(t, q)$. The series $F_{hv}(t, q) = F_{\geq \langle h, v\rangle}$ and $F_{d_1d_2}(t, q) = F_{\geq \langle d_1, d_2\rangle}$ are also computed in the following sections.

%
%
\section{Rotational symmetry classes}

We now study the convex polyominoes admitting a rotational symmetry. Observe that if a polyomino is $r$-symmetric, where $r$ denotes a $90^\circ$ rotation, then it is also $r^3$-symmetric and conversely, since $r^3 = r^{-1}$. In other words, $\Fix_\mathcal{P}(r) = \Fix_\mathcal{P}(r^3)$ and $F_r(t, q) = F_{r^3}(t, q)$. Hence it will be sufficient to determine the two generating series $F_r(t, q)$ and $F_{r^2}(t, q)$ of $r$-symmetric and of $r^2$-symmetric convex polyominoes, according to half-perimeter and area.

%
%
\subsection{$r$-Symmetric convex polyominoes}

\begin{figure}[!ht]
\begin{center}
\includegraphics{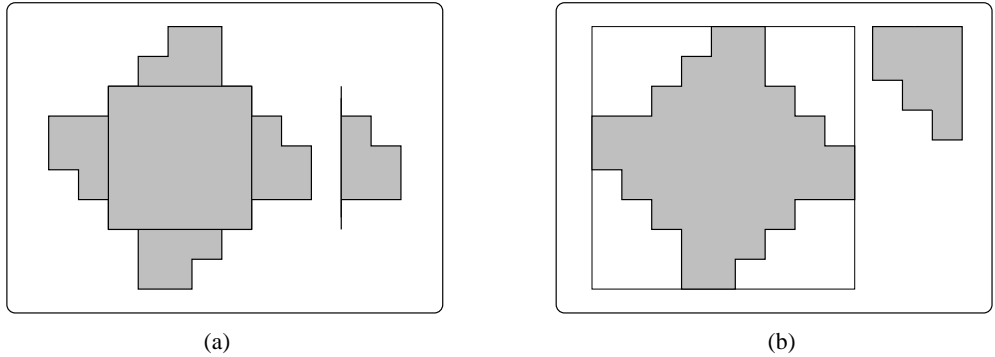}
\caption{\small Two ways of constructing a convex $r$-symmetric polyomino.}
\label{rSymDef}
\end{center}
\end{figure}

Observe that the minimal bounding rectangle of an $r$-symmetric convex polyomino is a square (see Figure~\ref{rSymDef}(b)). In particular, the half-perimeter is always even. There is also a maximal square fitting inside the polyomino (see Figure~\ref{rSymDef}(a)).

A first way to construct an $r$-symmetric convex polyomino $P$ is to start with its maximal fitting square $S$ and to glue the same stack polyomino $T$, after appropriate rotations, on each of its sides. The added stack can have empty parts on each side. The generating series $T_0(x,y,q)$ of such stacks is given in Subsection~\ref{sec:stack} (see formulas~(\ref{eqn:stack1}), (\ref{eqn:stack2}) and (\ref{eqn:stacks})).

\begin{prop}
\label{SerieGenRSym1}
The generating function for convex $r$-symmetric polyominoes \(F_r(t,q)\) is given by
\begin{equation}
\label{SerieR1}
F_r(t, q) = T_0(t^4,\  (y^m \mapsto t^{2m}q^{m^2})_{\, m \geq 1}\,,\  q^4)\,,
\end{equation}
\end{prop}

\begin{proof}
Suppose that the stack $T$ has height $m$, width $n$ and area $a$. Then the polyomino $P$ will have half-perimeter $2m+4n$ and area $m^2+4a$. The ``substitution'' (\ref{SerieR1}) gives exactly this.
\end{proof}

\begin{figure}[!ht]
\begin{center}
\includegraphics{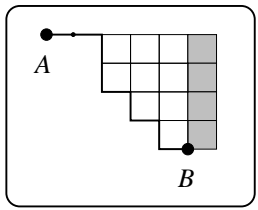}
\caption{\small Extended Ferrers diagram.}
\label{Fpasm}
\end{center}
\end{figure}

A second way to construct an $r$-symmetric convex polyomino $P$, of half-perimeter $2m$, say, is to start with the minimal bounding square, of side $m$, and to remove from each corner a given (appropriately rotated) Ferrers diagram $F$ with half-perimeter at most $m-1$ (see Figure~\ref{rSymDef}(b)). Denote by $D_m(q)$ the area generating series of Ferrers diagrams with half-perimeter less than or equal to $m$. If the perimeter of a diagram is less than $m$, then we can force it to be $m$ by adding empty columns to the diagram (see Figure~\ref{Fpasm}). We will call such diagrams, \textsl{extended} Ferrers diagrams. Thus $D_m(q)$ is the area generating function of the extended Ferrers diagrams with half-perimeter $m$. By convention $D_0(q) = 1$, and $D_1(q) = 1$ (one empty column).

\begin{prop}
\label{FormeCloseD}
For $m \geq 1$,
\begin{equation}
\label{DiagFm}
D_m(q) = \sum_{i=0}^{m-1} q^i  \left[ {m-1}\atop{i}\right]_q\,.
\end{equation}
\end{prop}

\begin{proof}
Let $i$ be the height of the first column on the right. Then this first column contributes for $q^i$ to the area. Removing that column yields a binomial path from $A$ to $B$, of length $m-1$ with $i$ vertical steps (see Figure~\ref{Fpasm}). The area generating function for such paths is given by the $q$-binomial coefficient $\left[{m-1}\atop{i}\right]_q$. Since $i$ can vary from $0$ to $m-1$, the conclusion follows.
\end{proof}

\begin{prop}
\label{propFerrers}
The polynomials $D_m(q)$ satisfy the recurrences
\begin{eqnarray}
\label{recurrenceF}
\mathbf{(a)}\quad D_m(q) & = & D_{m-1}(q) + q^{m-1} \sum_{h=0}^{m-2}D_h(q)\,, \quad m \geq 1\label{RRec1}\,,\\
\mathbf{(b)}\quad D_m(q) & = & (1+q)D_{m-1}(q) + (q^{m-1}-q)D_{m-2}(q)\,, \quad m \geq 2\,,\label{RRec2}
\end{eqnarray}
with $D_0(q) = 1 = D_1(q)$.
\end{prop}

\begin{figure}[!ht]
\begin{center}
\includegraphics{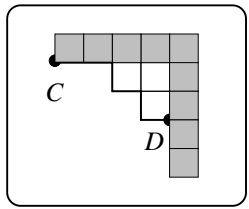}
\caption{\small Outer border and remaining path in a Ferrers diagram.}
\label{RecurF}
\end{center}
\end{figure}

\begin{proof}
\textbf{(a)}\ An extended Ferrers diagram with half-perimeter $m$ can take two forms:
\begin{enumerate}
\item
There is at least one empty column on the left end of the diagram. Truncating the first empty column gives $D_{m-1}(q)$.
\item
The diagram has no empty columns, i.e. is a normal Ferrers diagram. In that case we can remove the outer border (see Figure~\ref{RecurF}). There are $m-1$ cells in the border, contributing for $q^{m-1}$ to the area. What remains is an extended Ferrers diagram with half-perimeter $h$ (i.e. the length of the path from $C$ to $D$), $h$ varying from 0 (if the initial diagram consisted of a single column) to $m-2$ (if the initial diagram only had one row or no row of width $1$). This gives the second term of (\ref{RRec1}).

\textbf{(b)}\ The recurrence (\ref{RRec2}) is an immediate consequence of (\ref{RRec1}).
\end{enumerate}
\end{proof}

We can now compute the generating series $F_r(t,q)$ for the convex $r$-symmetric polyominoes.

\begin{prop}
The generating series \(F_r(t,q)\) for the $r$-symmetric convex polyominoes is given by
\begin{equation}
\label{SerieR2}
F_r(t,q) = \sum_{m \geq 1}t^{2m} q^{m^2}D_{m-1}(q^{-4})\,,
\end{equation}
where $D_{m-1}(q)$ is the generating series of the extended Ferrers diagrams with half-perimeter $m-1$. 
\end{prop}

\begin{proof}
This follows immediately from the second construction.
\end{proof}

If we denote by $a_{2m}(q)$ the generating series for the $r$-symmetric convex polyominoes with half-perimeter $2m$, then we get, for $m \geq 1$,
\begin{equation}
a_{2m}(q) = q^{m^2} D_{m-1}(q^{-4})\,.
\end{equation}
We have, for example, $a_0(q) = 1$ (by convention), $a_2(q) = q$ and $a_4(q) = q^4$. We then deduce from Proposition~\ref{FormeCloseD} the following closed form for these polynomials:

\begin{prop}
The area generating polynomial for $r$-symmetric convex polyominoes with half-perimeter $2m$ is given by
\begin{equation}
\label{FormeClosea}
a_{2m}(q) = q^{m^2} \left [ \sum_{i=0}^{m-2}q^i\qbinomial{m-2}{i}{q}\right]_{q \mapsto q^{-4}}\,,
\end{equation}
for $m \geq 2$.\hfill $\square$
\end{prop}

\begin{cor}
Let \(a_{2m} = a_{2m}(1)\), the number of convex $r$-symmetric polyominoes with half-perimeter $2m$. Then $a_2 = 1$ and 
\begin{equation}
a_{2m} = 2^{m-2}
\end{equation}
for $m \geq 2$\,.\hfill $\square$
\end{cor}

The following recurrences are useful in order to compute explicitly the polynomials $a_{2m}(q)$.

\begin{prop}
The polynomials $a_{2m}(q)$ satisfy the recurrences
\begin{eqnarray}
\mathbf{(a)}\quad a_{2m}(q) & = & q^{2m-1}a_{2m-2}(q) + \sum_{j=1}^{m-2}q^{(m-2)^2 + 4 -j^2} a_{2j}(q)\,,\quad m \geq 1\,,\label{RRec3}\\
\mathbf{(b)}\quad a_{2m}(q) & = & q^{2m-5}(1+q^4)a_{2m-2}(q) + q^4(1-q^{4m-12})a_{2m-4}(q)\,, \quad m \geq 3\,,\label{RRec4}
\end{eqnarray}
with initial conditions $a_0(q) = 1$, $a_2(q) = q$ and $a_4(q) = q^4$.
\end{prop}

\begin{proof}
This follows immediately from Proposition~\ref{propFerrers}.
\end{proof}

Here are the first few terms of the generating series \(F_r(t,q)\) for convex $r$-symmetric polyominoes:
\begin{equation}
F_r(t,q) = t^2q+t^4q^4+t^6q^5+t^6q^9+2\,t^8q^8+t^{10}q^9+t^8q^{12}+ \ldots
\end{equation}

%
%
\subsection{$r^2$-Symmetric convex polyominoes}

$r^2$-symmetric convex polyominoes fall into two categories: those with even width and those with odd width. Thus, if $F_{r^2}(x,y,q)$ denotes the generating series for $r^2$-symmetric convex polyominoes according to width, height and area, we have
\begin{equation}
\label{SerieRCarree}
F_{r^2}(x,y,q) = F_{r^2,\,e}(x,y,q) +F_{r^2,\,o}(x,y,q)\,,
\end{equation}
where $F_{r^2,\,e}$ (resp. \(F_{r^2,\,o}\)) is the generating series for $r^2$-symmetric convex polyominoes of even width (resp. odd width). We calculate these series in the next two propositions. Recall that $T(x,y,q)$ is the generating function of stack polyominoes and $D(s,x,y,q)$, that of directed (one source) convex polyominoes, where $s$ is an extra variable which marks the length of the leftmost column (see Subsection~\ref{sec:dirconv}).

\begin{figure}[ht]
\begin{center}
\includegraphics{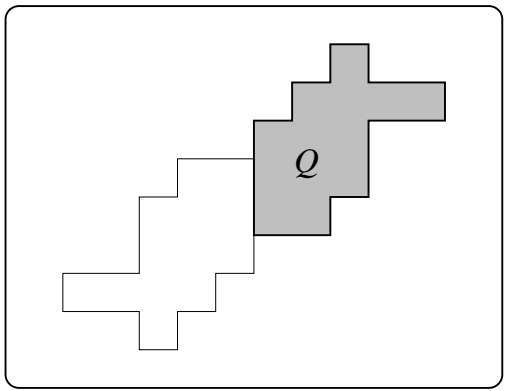}
\caption{\small An $r^2$-symmetric convex polyomino $P$ and its fundamental region $Q$.}
\label{r2sym}
\end{center}
\end{figure}

\begin{prop}
Let \(F_{r^2,\,e}(x,y,q)\) be the generating series for $r^2$-symmetric convex polyominoes of even width. Then we have
\begin{equation}
F_{r^2,\,e}(x,y,q) = \frac{2}{1-y} \left( D(\frac{1}{y}, x^2,y^2, q^2) - D(1, x^2, y^2, q^2) \right) - T(x^2, y, q^2)\,.
\end{equation}
\end{prop}

\begin{figure}[!ht]
\begin{center}
\includegraphics{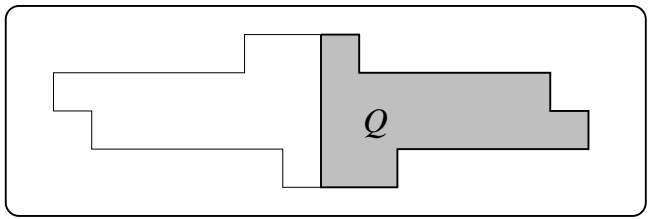}
\caption{\small Convex $r^2$-symmetric polyomino whose fundamental region is a stack.}
\label{R2-symTas}
\end{center}
\end{figure}

\begin{proof}
Let $P$ be an $r^2$-symmetric convex polyomino of even width. We define the fundamental region of $P$ to be the right half of $P$. Call this polyomino $Q$ (see Figure~\ref{r2sym}). Three cases can occur according to the privileged direction taken by $Q$ with respect to $P$.

In the first case, the fundamental region $Q$ is directed North-East (Figure~\ref{r2sym}). To recover $P$, we rotate a copy of $Q$ by $180^\circ$ and glue the result $\bar{Q}$ to $Q$ along the leftmost column. If this column has length equal to $k$, there will be $k$ possible positions for $\bar{Q}$ relative to $Q$. The substitution $s \mapsto 1/y$, $x \mapsto x^2$, $y \mapsto y^2$ and $q \mapsto q^2$ in the generating series $D(s,x,y,q)$ of directed convex polyominoes corresponds to the highest position of $\bar{Q}$, which minimizes the overall height of $P$. All the possible positions will be accounted for by multiplying by $(1+y+\ldots + y^{k-1})$. In other words, the substitution to make in $D(s,x^2,y^2,q^2)$ is
\begin{equation}
s^k \mapsto \frac{1 + y + \ldots + y^{k-1}}{y^k} = \frac{1}{1-y}\left(\frac{1}{y^k}-1\right)\,.
\end{equation}

Summing over all possible $k$'s, we find that the generating series for $r^2$-symmetric convex polyominoes in the first class is
\begin{equation}
\label{eqnR2}
\frac{1}{1-y} \left(D(\frac{1}{y},x^2,y^2,q^2) - D(1, x^2, y^2, q^2) \right)
\end{equation}

In the second case, the fundamental region $Q$ is directed South-East. Using an horizontal reflection, we see that the $r^2$-symmetric convex polyominoes in this class are also enumerated by (\ref{eqnR2}).

In the third case, the fundamental region $Q$ does not take a definite North or South direction. $Q$ is then a stack polyomino, and the rotated image $\bar{Q}$ is glued to $Q$ in the minimum height position (see Figure~\ref{R2-symTas}). Hence these $r^2$-symmetric convex polyominoes are enumerated by $T(x^2,y,q^2)$. However, this class is already contained both in the first and the second classes, so that this term must in fact be subtracted from the others.
\end{proof}

\begin{prop}
The generating series $F_{r^2,\,o}(x,y,q)$ for convex $r^2$-symmetric polyominoes with odd width is given by
\begin{equation}
F_{r^2,\,o}(x,y,q) = \frac{2}{x} D(\frac{1}{yq}, x^2, y^2, q^2) - \frac{1}{x} T(x^2, \frac{y}{q}, q^2)\,.
\end{equation}
\end{prop}

\begin{figure}[!ht]
\begin{center}
\includegraphics{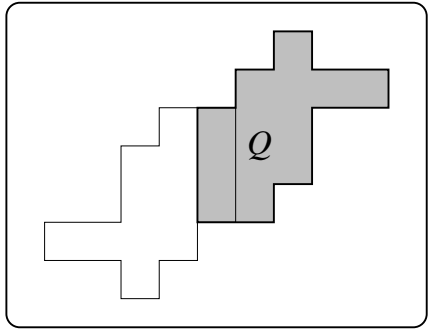}
\caption{\small Convex $r^2$-symmetric polyomino with odd width.}
\label{R2SymImp}
\end{center}
\end{figure}

\begin{proof}
The proof is similar to the previous one. The main difference is that only one glueing position of $\bar{Q}$ to $Q$ is admissible and that furthermore the leftmost column of $Q$ and its rotated image in $\bar{Q}$ are superimposed to yield an odd width (see Figure~\ref{R2SymImp}). Details are left to the reader.
\end{proof}

We give here the first few terms of the generating series $F_{r^2}(t,t,q)$ according to half-perimeter and area for $r^2$-symmetric convex polyominoes:
\begin{equation}
F_{r^2}(t,t,q) = t^{2}q+2\,t^{3}q^{2}+2\,t^{4}q^{3}+t^{4}q^{4}+6\,t^{5
}q^{4}+2\,t^{5}q^{6}+7\,t^{6}q^{5}+ \ldots
\end{equation}

%
%
\subsection{Rotation-type convex polyominoes}

We can now apply formula~(\ref{BurnsideL}) to the case where $\mathcal{P}$ is the class $C$ of convex (translation-type) polyominoes and $G$ is the subgroup $\mathfrak{C}_4$ of $\mathfrak{D}_4$, of rotations, since rotation-type polyominoes are precisely defined as orbits of polyominoes relative to $\mathfrak{C}_4$. We thus get
\begin{equation}
|C/\mathfrak{C}_4| = \frac{1}{4}\left( |\Fix(1)| + |\Fix(r)| + |\Fix(r^2)| + |\Fix(r^3)|
\right)
\end{equation}

Note that $|\Fix(1)|$ is the generating series $C(t,t,q)$ of convex polyominoes (see Subsection~\ref{sec:dirconv}). $|\Fix(r)|_w = |\Fix(r^3)|_w = F_r(t,q)$ is given by equation~(\ref{SerieR1}) or equation~(\ref{SerieR2}) and $|\Fix(r^2)|_w = F_{r^2}(t,t,q)$ is given by equation~(\ref{SerieRCarree}).

\begin{prop}The half-perimeter and area generating series for convex rotation-type polyominoes $(C/\mathfrak{C}_4)(t,q)$ is given by
\begin{equation}
(C/\mathfrak{C}_4)(t,q) = \frac{1}{4} \left( C(t,t,q) + 2F_r(t,q) + F_{r^2}(t,t,q) \right)\,.
\end{equation}\hfill $\square$
\end{prop}

The first few terms of $(C/\mathfrak{C}_4)(t,q)$ are given in Appendix B.

%
%
\section{Horizontal and vertical symmetry classes}

%
%
\subsection{Horizontal or vertical symmetry}

We first note that there is an obvious bijection between the horizontally and vertically symmetric polyominoes: rotation by a right angle. We can thus restrict our study to the $v$-symmetric polyominoes (those left unchanged by a reflection about some vertical axis). 

Define the fundamental region of a $v$-symmetric convex polyomino $P$ to be the part $Q$ of $P$ on the right hand side of the symmetry axis (including the central column in the odd width case). This region is a stack polyomino (see Figure~\ref{VertPetit}).

\begin{prop}
\label{SerieGenVSym}
Let $F_v(t,q)$ be the generating series for the $v$-symmetric convex polyominoes. Then
\begin{equation}
F_v(t,q) = T(t^2, t, q^2) + \frac{1}{t}T(t^2, \frac{t}{q}, q^2)\,,
\end{equation}
where $T(x,y,q)$ is the generating series for stack polyominoes.
\end{prop}

\begin{proof}
Let $P$ be a convex $v$-symmetric polyomino. There are two possibilities:
\begin{enumerate}
\item
\label{Etape1}
The width of $P$ is even (see Figure~\ref{VertPetit}(a)), in which case the symmetry axis lies between two columns of $P$. Consider the stack polyomino $Q$ right of the axis. Reflecting it about the axis doubles the width and the area, but the height remains unchanged. Hence the generating series in this case is given by $T(t^2, t, q^2)$.

\begin{figure}[!ht]
\begin{center}
\includegraphics{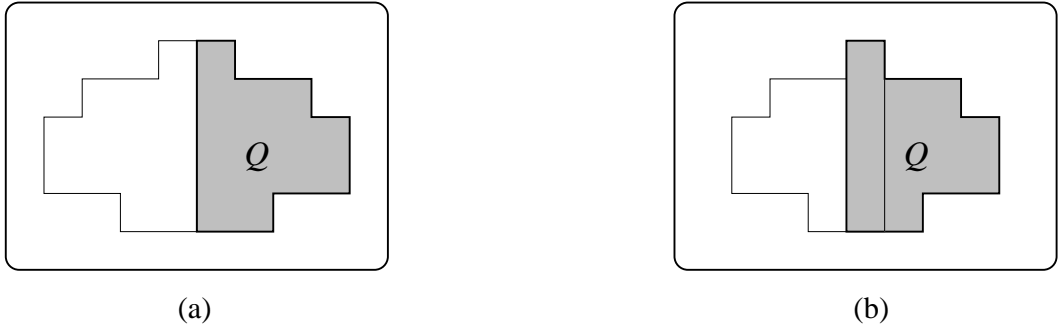}
\caption{\small Convex $v$-symmetric polyomino with: (a) even width, (b) odd width.}
\label{VertPetit}
\end{center}
\end{figure}

\item
The width of $P$ is odd (see Figure~\ref{VertPetit}(b)), in which case the symmetry axis passes through the central column. To construct $P$, reflect the fundamental region and glue the two parts together, taking care to superimpose the two central columns (identical by symmetry) to get an odd-width polyomino. The correct generating series for this case is then $\frac{1}{t}T(t^2, \frac{t}{q}, q^2)$.

\end{enumerate}
\end{proof}

Here are the first few terms of the series $F_v(t,q)$:
\begin{equation}
F_v(t,q) = t^{2}q+2\,t^{3} q^{2}+2\,t^{4}q^{3}+t^{4}q^{4}+4\,t^{5}q^{4}+2\,t^{5 }q^{6}+5\,t^{6}q^{5}+
\ldots
\end{equation}

%
%
\subsection{Horizontal and vertical symmetry}

Let us examine now the set of convex polyominoes that are both horizontally and vertically symmetric. Such a polyomino $P$ will be called a $hv$-symmetric polyomino. Define the fundamental region $Q$ of $P$ to be the bottom-right part of $P$, that is the part right of the vertical symmetry axis and down the horizontal symmetry axis. The region $Q$ is a partition polyomino (see Figure~\ref{VHpp} and Subsection~\ref{sec:partition}). 

\begin{prop}
\label{SerieGenVHSym}
The generating series $F_{hv}(t,q)$ for the convex $hv$-symmetric polyominoes is given by
\begin{equation}
\label{VHseries}
F_{hv}(t,q) = P(t^2, t^2, q^4) + 2 \frac{1}{t} P(\frac{t^2}{q^2}, t^2, q^4) + \frac{q}{t^2}
P(\frac{t^2}{q^2}, \frac{t^2}{q^2}, q^4)\,,
\end{equation}
where $P(x,y,q)$ is the generating series for partition polyominoes.
\end{prop}

\begin{proof}
Let $P$ be a convex $hv$-symmetric polyomino. $P$ falls into one of four disjoint cases, according to the parities of its height and width.

\begin{enumerate}
\item
$P$ has even height and even width. Then we simply have to glue together four copies of the fundamental region $Q$ (see Figure~\ref{VHpp}).
\begin{figure}[!ht]
\begin{center}
\includegraphics{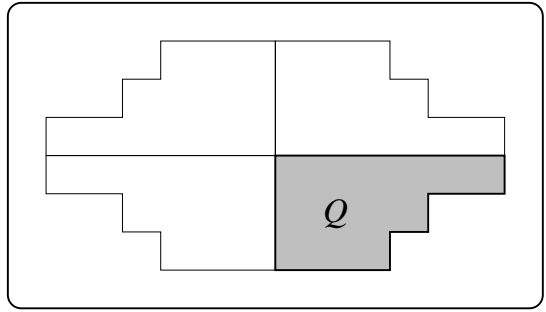}
\caption{\small Convex $hv$-symmetric polyomino with even width and even height.}
\label{VHpp}
\end{center}
\end{figure}
Replacing $x$ and $y$ by $t^2$ we will get the correct perimeter value for the whole object, and replacing $q$ by $q^4$ will give the right area. So the generating series $F_{hv,\,ee}(t,q)$ for this type of polyominoes is
\begin{equation}
\label{vhpp}
F_{hv,\,ee}(t,q) = P(t^2, t^2, q^4)\,.
\end{equation}

\item
$P$ has even width and odd height. Then the horizontal symmetry axis passes through the middle row (see Figure~\ref{VHpi}(a)). We first duplicate four times the fundamental region, glue them together as above, but merge the two central rows. The correct adjustment requires replacing $x$ by $t^2/q^2$ and dividing by $t$. The generating series $F_{hv,\,eo}(t,q)$ for these polyominoes is thus
\begin{equation}
\label{vhpi}
F_{hv,\,eo}(t,q) = \frac{1}{t} P(\frac{t^2}{q^2}, t^2, q^4)\,.
\end{equation}

\begin{figure}[!ht]
\begin{center}
\includegraphics{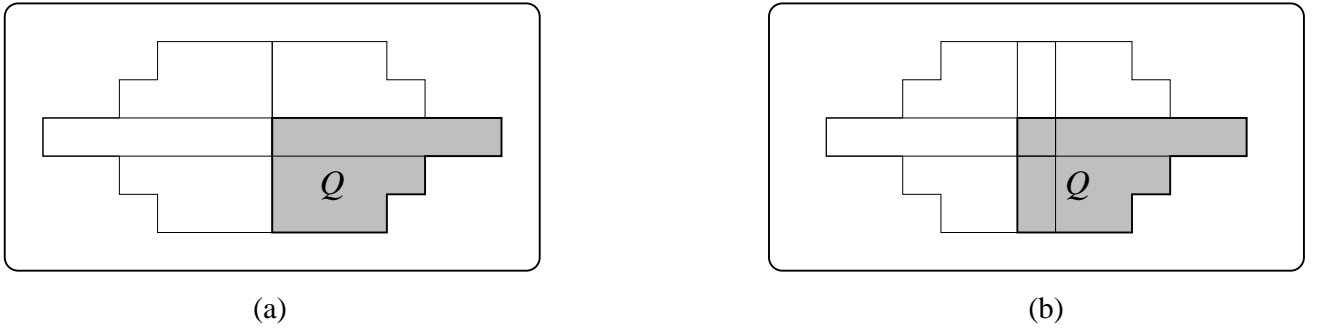}
\caption{\small Convex $hv$-symmetric polyomino with odd height and: (a) even width, (b) odd width.}
\label{VHpi}
\end{center}
\end{figure}

\item
$P$ has odd width and even height. Using a $90^\circ$ rotation, we see that this case reduces to the preceeding one by interchanging $x$ and $y$ in $P(x,y,q)$. Since $P(x,y,q)$ is symmetric in $x$ and $y$, the generating series $F_{hv,\,oe}(t,q)$ for this case is also given by
\begin{equation}
\label{vhip}
F_{hv,\,oe}(t,q) = \frac{1}{t} P(\frac{t^2}{q^2}, t^2, q^4)\,.
\end{equation}

\item
$P$ has odd width and odd height (see Figure~\ref{VHpi}(b)). The generating series $F_{hv,\,oo}(t,q)$ for this class of polyominoes is given by 
\begin{equation}
\label{vhii}
F_{hv,\,oo}(t,q) = \frac{q}{t^2} P(\frac{t^2}{q^2}, \frac{t^2}{q^2}, q^4)\,.
\end{equation}
Summing the four series yields~(\ref{VHseries}).
\end{enumerate}
\end{proof}

Let $f_n^{ee}(q)$, $f_n^{eo}(q)$, $f_n^{oe}(q)$ and $f_n^{oo}(q)$ denote the area generating polynomials for these classes of $hv$-symmetric polyominoes with half-perimeter equal to $n$. In other words,
\begin{displaymath}\left.
\begin{array}{ll}
F_{hv, ee}(t,q) = \sum_{n \geq 1}f_n^{ee}(q)\,t^n\,, & F_{hv, eo}(t,q) = \sum_{n \geq 1}f_n^{eo}(q)\,t^n\,,\\
F_{hv, oe}(t,q) = \sum_{n \geq 1}f_n^{oe}(q)\,t^n\,, & F_{hv, oo}(t,q) = \sum_{n \geq 1}f_n^{oo}(q)\,t^n\,.
\end{array}\right.
\end{displaymath} 

Obviously, we have, for $n \geq 1$,
\begin{eqnarray}
f_{n+2}^{ee}(q) & = & q^{2n-1}f_n^{oo}(q)\,, \label{fee}\\
f_{n}^{oe}(q) & = & f_n^{eo}(q)\,, \label{feo}
\end{eqnarray}
so that we can concentrate our attention on the polynomials $f_{n}^{oo}(q)$ and $f_{n}^{eo}(q)$.

\begin{prop}
The area generating polynomials $f_{n}^{oo}(q)$ and $f_{n}^{eo}(q)$ are given by
\begin{equation}
\label{eqn:oopoly}
f_{n}^{oo}(q) = \left\{ \begin{array}{ll}
q^{n-1}\displaystyle \sum_{i=0}^{(n-2)/2}\qbinomial{(n-2)/2}{i}{{q^4}} & \textrm{if $n$ is even}\,,\\
0 & \textrm{if $n$ is odd}\,,
\end{array}\right.
\end{equation}

\begin{equation}
\label{eqn:eopoly}
f_{n}^{eo}(q) = \left\{ \begin{array}{ll}
\displaystyle q^{n-1}\sum_{i=0}^{(n-3)/2}q^{2i}\qbinomial{(n-3)/2}{i}{{q^4}} & \textrm{if $n$ is odd}\,,\\
0 & \textrm{if $n$ is even}\,.
\end{array}\right.
\end{equation}

Moreover they satisfy the following recurrences: for $n \geq 5$,
\begin{eqnarray}
f_{n}^{oo}(q) & = & 2q^2\,f_{n-2}^{oo}(q) - q^4(1-q^{2n-8})\,f_{n-4}^{oo}(q)\,, \label{eqn:oorec}\\
&  & \nonumber \\
f_{n}^{eo}(q) & = & q^2(1+q^2)\,f_{n-2}^{eo}(q) - q^6(1-q^{2n-10})\,f_{n-4}^{eo}\,, \label{eqn:eorec}
\end{eqnarray}

with initial conditions
\begin{displaymath}\left.
\begin{array}{llll}
f_{1}^{oo}(q) = 0\,, & f_{2}^{oo}(q) = q\,, & f_{3}^{oo}(q) = 0\,, & f_{4}^{oo}(q) = 2q^3\,, \\
f_{1}^{eo}(q) = 0\,, & f_{2}^{eo}(q) = 0\,, & f_{3}^{eo}(q) = q^2\,, & f_{4}^{eo}(q) = 0\,.
\end{array}\right.
\end{displaymath}
\end{prop}

\begin{proof}The explicit formulas~(\ref{eqn:oopoly}) and (\ref{eqn:eopoly}) can be easily derived from the generating series or from a direct combinatorial argument. The recurrences~(\ref{eqn:oorec}) and (\ref{eqn:eorec}) can be checked by substituting the explicit formulas.
\end{proof}

\begin{cor}The number $f_n$ of $hv$-symmetric convex polyominoes with half-perimeter $n$ is $0$ if $n = 1$, $1$ if $n = 2$, and for $n \geq 3$, is given by
\begin{eqnarray}
f_n = f_n^{oo}(1) + f_n^{ee}(1) & = & 3\cdot 2^{(n/2)-2}  \quad \textrm{for $n$ even}\,, \label{eqn:numeven}\\
f_n = f_n^{eo}(1) + f_n^{oe}(1) & = & 2^{(n-1)/2}  \quad \ \,\textrm{for $n$ odd}\,. \label{eqn:numodd}
\end{eqnarray}

Furthermore, we have
\begin{equation}
F_{hv}(t,1) = \sum_{n \geq 1}f_n\,t^n = \frac{t(t+1)^2}{1-2t^2}\,.
\end{equation}\hfill $\square$
\end{cor}

We give here the first few terms of the series $F_{hv}(t,q)$: 
\begin{equation}
F_{hv}(t,q) = t^{2}q+2\,t^{3}q^{2}+2\,t^{4}q^{3}+t^{4}q^{4 }+2\,t^{5}q^{4}+2\,t^{5}q^{6}+3\,t^{6}q^{5} +
\ldots
\end{equation}

%
%
\section{Diagonal symmetry classes}
\label{sec:diagsymclass}

%
%
\subsection{Diagonally symmetric convex polyominoes}

\begin{figure}[ht]
\begin{center}
\includegraphics{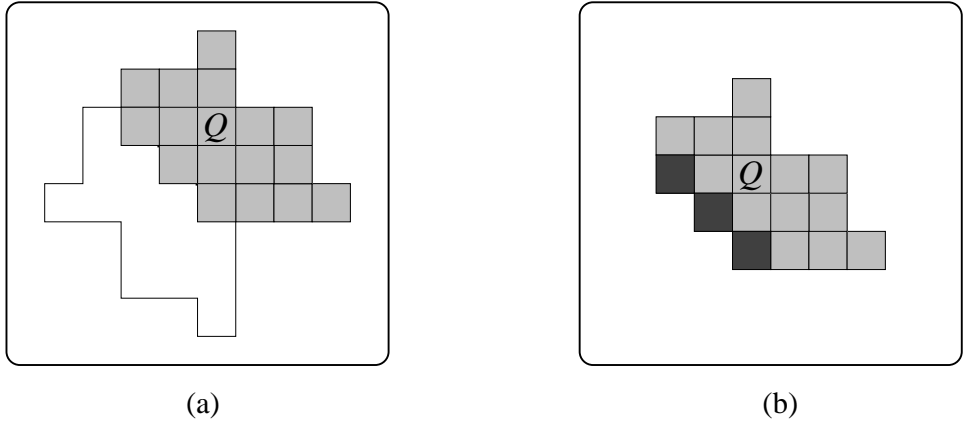}
\caption{\small Diagonally symmetric convex polyomino $P$ and its fundamental region $Q$.}
\label{DSym}
\end{center}
\end{figure}

Our next goal is to compute the generating series $F_d(t,q)$ of $d$-symmetric convex polyominoes. There is a simple bijection between $d_1$-symmetric and $d_2$-symmetric polyominoes, given by the $90^\circ$ rotation $r$, so that $F_{d_1}(t,q) = F_{d_2}(t,q)$. We choose $d = d_1$, the diagonal along the line $y=-x$ and define the \textsl{fundamental region} $Q$ of a $d$-symmetric convex polyomino $P$ as the part of $P$ located on the North-East side of the diagonal symmetry line, including the cells on the diagonal (see Figure~\ref{DSym}).

Observe that the fundamental region $Q$ of a $d$-symmetric convex polyomino $P$ is a (diagonally) shifted directed (diagonal source) convex polyomino, called \textsl{$D_S$-polyomino} for short. These should not be confused with the neighboring class of directed diagonally convex polyominoes, also known as fully directed compact lattice animals (see \cite{VPNMS2} and \cite{VPNMS}), and a different enumerating method has to be used.

We introduce an extra variable $z$ to mark the length of the diagonal source. For example, the $D_S$-polyomino in Figure~\ref{DSym}(b) has global weight $x^6y^5z^3q^{17}$. Let $D_S(x,y,z,q)$ denote the generating series of $D_S$-polyominoes. Before computing this series, we can state the following:

\begin{prop}The half-perimeter and area generating series $F_d(t,q)$ of $d$-symmetric convex polyominoes is given by
\begin{equation}
\label{eqn:diagsubs}
F_d(t,q) = D_S(t^2, t^2, \frac{1}{t^2q}, q^2)\,,
\end{equation}
where $D_S(x,y,z,q)$ is the generating series of shifted directed convex polyominoes.
\end{prop}

\begin{proof}Observe that if the fundamental region $Q$ of  a $d$-symmetric polyomino $P$ has width $m$, height $n$, diagonal source length $k$ and area $a$, then $P$ has half-perimeter equal to $2m+2n-2k$ and area equal to $2a-k$. This is achieved precisely by the substitution $x = t^2$, $y = t^2$, $z = 1/t^2q$ and $q = q^2$ as in formula~(\ref{eqn:diagsubs}).
\end{proof}

Here are the first few terms:
\begin{equation}
F_{d}(t,q) = t^{2}q + 2\,t^{4}q^{3} + t^4q^4 + 5\,t^6q^5 + 4\,t^6q^6 + 2\,t^6q^7 + 2\,t^6q^8 + \ldots
\end{equation}

%
%
\subsection{Shifted directed convex ($D_S$) polyominoes}

We now outline some basic enumerative results concerning the class $D_S$ of shifted directed convex polyominoes. This class is studied in more detail in a separate paper \cite{AJGPL}. 

\begin{figure}[ht]
\begin{center}
\includegraphics{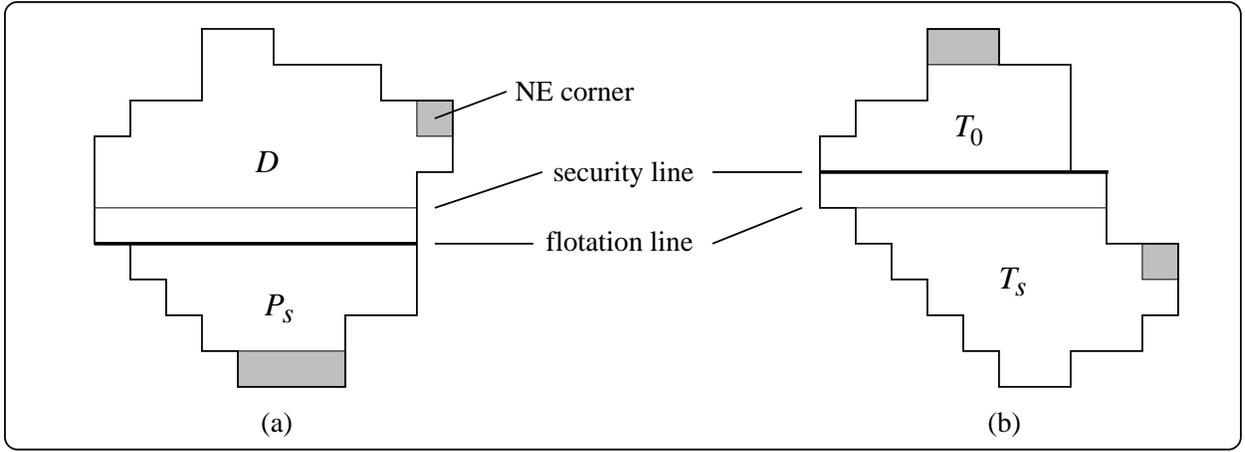}
\caption{\small $D_S$-polyominoes of types $Y_1$ (case (a)) and $Y_2$ (case (b)).}
\label{figy1y2}
\end{center}
\end{figure}

Let $Q$ be a $D_S$-polyomino. Imagining $Q$ as a kind of floating device, we define the \textsl{flotation row} of $Q$ as the (horizontal) row containing the top diagonal source cell. The two bounding lines, below and above the flotation row are called respectively the \textsl{flotation line} and the \textsl{security line} of $Q$. We also define the \textsl{North-East corner} of $Q$ as the topmost cell of the rightmost column of $Q$ (see Figure~\ref{figy1y2}).

We now partition the class of $D_S$-polyominoes into two subclasses denoted by $Y_1$ and $Y_2$, according to the relative position of the North-East corner with respect to the flotation line: above for $Y_1$ and below for $Y_2$. We also denote their generating series by the same names and we have
\begin{equation}
D_S(x,y,z,q) = Y_1(x,y,z,q) + Y_2(x,y,z,q)\,.
\end{equation}

In order to enumerate $Y_1$-polyominoes, we introduce a new variable $v$ to mark the width of the bottom row. For example, the polyomino $Q$ of Figure~\ref{DSym}(b) has bottom width equal to $4$. Let $Y_1(v,x,y,z,q)$ denote the corresponding generating function. Observe that the (one source) directed convex polyominoes are a special case. The corresponding generating function is obtained from $D(s,x,y,q)$ (see Subsection~\ref{sec:dirconv} and \cite{MBMa}) by setting $s = v$, applying a diagonal reflection. Extracting the coefficient of $z$ in $Y_1$ gives
\begin{equation}
D(v,x,y,q) = [z]\,Y_1(v,x,y,z,q)\,.
\end{equation}

Now any $D_S$-polyomino $Q$ in the class $Y_1$ can be constructed by glueing together along the flotation line a directed convex polyomino $D$, of bottom width $m \geq 1$, with a (possibly empty) diagonally decreasing stack of rows, i.e. a shifted partition $P_S$ with parts of length less than $m$ (see Figure~\ref{figy1y2}(a)). This gives the following:

\begin{prop}Let $D_{m}(x,y,q) = [v^m]\,D(v,x,y,q)$ be the generating series for directed convex polyominoes of bottom width $m$. We then have
\begin{equation}
Y_1(x,y,z,q) = \sum_{m \geq 1}zD_{m}(x,y,q)(1+yzq)\cdots (1+yzq^{m-1}).
\end{equation}\hfill $\square$
\end{prop}

It is possible to get a more closed form for this generating function, using the following functional equation:
\begin{prop}
The generating function $Y_1(v) = Y_1(v,x,y,z,q)$ of $Y_1$-polyominoes is characterized by the equation
\begin{equation}
\label{eqn:y1eq}
Y_1(v) = zD(v,x,y,q) + \frac{yz}{1-vq}(vqY_1(1) - Y_1(vq))\,.
\end{equation}
\end{prop}

\begin{proof}This equation says that any $Y_1$-polyomino can be obtained as follows: start with a directed convex polyomino $D$ and then successively glue, on the bottom, rows that follow the diagonal on the left hand side but do not exceed the previous row on the right hand side. The term $yzvq(1-vq)^{-1}Y_1(1)$ accounts for the addition of a row of arbitrary length whereas $yz(1-vq)^{-1}Y_1(vq)$ enumerates the cases where the added row does in fact exceed on the right hand side, and thus should be subtracted.
\end{proof}

\begin{cor}The generating series for $Y_1$-polyominoes is given by
\begin{equation}
Y_1(x,y,z,q) = \frac{\displaystyle \sum_{n \geq 0}zD(q^n,x,y,q)\frac{(-yz)^n}{(q)_n}}{\displaystyle 1 + \sum_{n \geq 1}q^{n-1}\frac{(-yz)^n}{(q)_n}}\,.
\end{equation}
\end{cor}

\begin{proof}
This is the solution $Y_1(1)$ of equation~(\ref{eqn:y1eq}) obtained by iteration. See Lemma 2.3 of \cite{MBMa}.
\end{proof}

We now study the subclass $Y_2$ of $D_S$-polyominoes. Let us introduce a variable $u$ to mark the top row width and let $Y_2(u,x,y,z,q)$ be the corresponding generating series. Cutting a $Y_2$-polyomino $P$ along the security line yields a stack polyomino $T$ on top and a shifted stack polyomino $T_S$ (see Subsection~\ref{sec:diagrowstack}) below the security line (see Figure~\ref{figy1y2}(b)). Note that $T_S$ is not a general diagonal row stack polyomino since the North-East corner is not allowed to be on the top row. Call $R$ the subclass of $T_S$ satisfying this restriction and observe that $T_S$-polyominoes whose North-East corner is on the top row are precisely shifted partitions (see Subsection~\ref{sec:partition}). We then have
\begin{equation}
R(u,x,y,q) = T_S(u,x,y,q) - P_S(xu,y,q)\,.
\end{equation}

Let $R_{m}(x,y,q)$ be the generating series of $R$-polyominoes having top width equal to $m$, i.e. set
\begin{equation}
R(u,x,y,q) = \sum_{m \geq 1}R_{m}(x,y,q)u^m.
\end{equation}

Also denote by $T_{0,m}(y,q)$ the generating series of stack polyominoes of width $m$ and having possibly empty columns (see Subsection~\ref{sec:stack}\,; note that these stacks have been rotated by $90^\circ$). We then have

\begin{prop}The generating series for $Y_2$-polyominoes is given by
\begin{equation}
Y_2(x,y,z,q) = \sum_{m \geq 1}R_{m}(x,yz,q)T_{0,m}(y,q)\,.
\end{equation}\hfill $\square$
\end{prop}

An alternate expression can be obtained for $Y_2$ by solving the following functional equation for $Y_2(u,x,y,z,q)$, with the help of Lemma 2.2 of \cite{MBMa}.

\begin{prop}
The series $Y_2(u) = Y_2(u,x,y,z,q)$ satisfies the functional equation
\begin{equation}
Y_2(u) = R(u,x,yz,q) + \frac{yuq}{1-uq}Y_2'(1) - \frac{yu^2q^2}{(1-uq)^2}(Y_2(1) - Y_2(uq))\,,
\end{equation}
where the derivative $Y_2'(u)$ is taken with respect to $u$.
\end{prop}

\begin{proof}
This functional equation states in fact that a $Y_2$-polyomino can be obtained as follows: start with a $R$-polyomino and then glue successively rows to the top. Each added row should not exceed the previous row on either side. This corresponds to the construction C of \cite{MBMa} where the general functional equation is stated and solved.
\end{proof}

Here are the first few terms of $D_S$:
\begin{equation}
D_S(x,y,z,q) = xyzq + x^2yzq^2 + xy^2zq^2 + x^3yzq^3 + x^2y^2z^2q^3 + 3\,x^2y^2zq^3 + xy^3zq^3 + \ldots 
\end{equation}

%
%
\subsection{Symmetry with respect to both diagonals}

A polyomino which is symmetric with respect to both diagonals is called \textsl{$d_1d_2$-symmetric}. Let $P$ be a $d_1d_2$-symmetric convex polyomino. Observe that the height of $P$ must be equal to its width. Two cases occur according to the parity of the width. Let $F_{d_1d_2,e}(t,q)$ and $F_{d_1d_2,o}(t,q)$ denote the half-perimeter and area generating functions for even and odd width $d_1d_2$-symmetric convex polyominoes respectively. We then have
\begin{equation}
\label{eqn:d1d2}
F_{d_1d_2}(t,q) = F_{d_1d_2,e}(t,q) +F_{d_1d_2,o}(t,q).
\end{equation}

\begin{figure}[!ht]
\begin{center}
\includegraphics{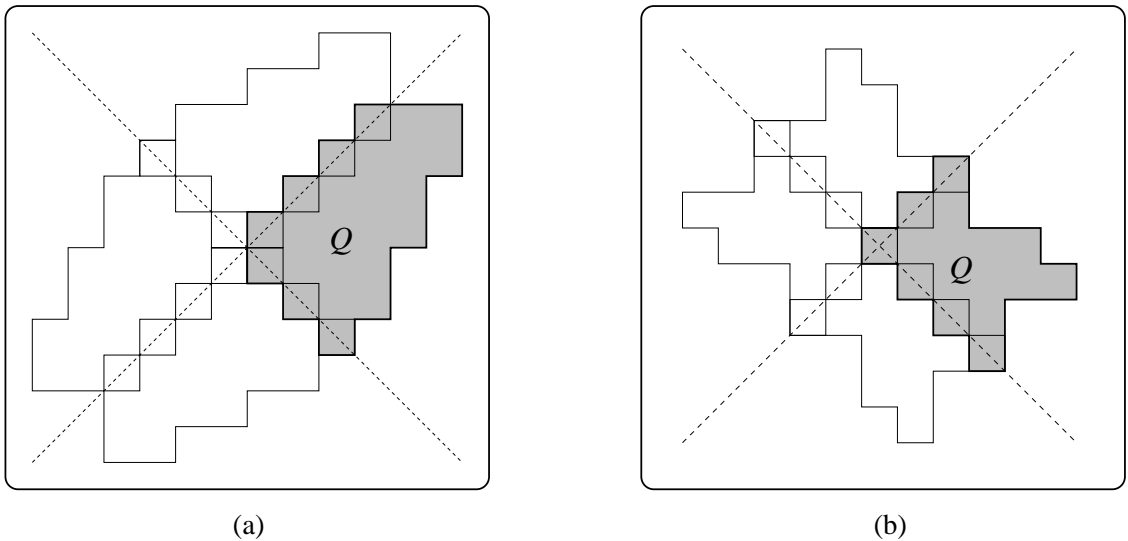}
\caption{\small Convex $d_1d_2$-symmetric polyomino with: (a) even width; (b) odd width.}
\label{DiagDiag}
\end{center}
\end{figure}

We define the \textsl{fundamental region} $Q$ \textsl{of} $P$ to be the part of $P$ which is bounded by the two diagonals and extends in the East direction, including the bounding diagonal cells (see Figure~\ref{DiagDiag}). Observe that $Q$ is a stack of rows with a double diagonal shift. We call these polyominoes \textsl{doubly shifted stacks ($T_{SS}$)} and qualify them as \textsl{even} in the even width case (Figure~\ref{DiagDiag}(a)) and \textsl{acute} in the odd case (Figure~\ref{DiagDiag}(b)). We use the variables $z$ to mark the length of the $d_1$-diagonal and $w$ for the $d_2$-diagonal, and introduce the generating series $E(x,z,w,q)$ and $A(x,z,w,q)$ of even and acute (respectively) doubly shifted stacks. Note that in the acute case, the extreme West cell contributes to both diagonal lengths. Before determining these series, we can state the following:

\begin{prop}The generating series for even and odd $d_1d_2$-symmetric convex polyominoes are given by
\begin{eqnarray}
F_{d_1d_2,e}(t,q) & = & E(t^4, \frac{1}{q^2}, \frac{1}{q^2}, q^4)\,, \label{eqn:E}\\
F_{d_1d_2,o}(t,q) & = & \frac{q}{t^2}A(t^4, \frac{1}{q^2}, \frac{1}{q^2}, q^4)\,. \label{eqn:A}
\end{eqnarray}
\end{prop}

\begin{proof}Suppose that the doubly shifted stack $Q$ has width $m$, $d_1$ and $d_2$ diagonal lengths $n_1$ and $n_2$ respectively, and area $a$. Then, in the even case, $P$ will have half-perimeter $4m$ and area $4a-2n_1-2n_2$, while in the odd case, the half-perimeter will be $4m-2$ and the area, $4a-2n_1-2n_2+1$. Since the correspondence between $P$ and $Q$ is bijective, we have (\ref{eqn:E}) and (\ref{eqn:A}). 
\end{proof}

\begin{figure}[!ht]
\begin{center}
\includegraphics{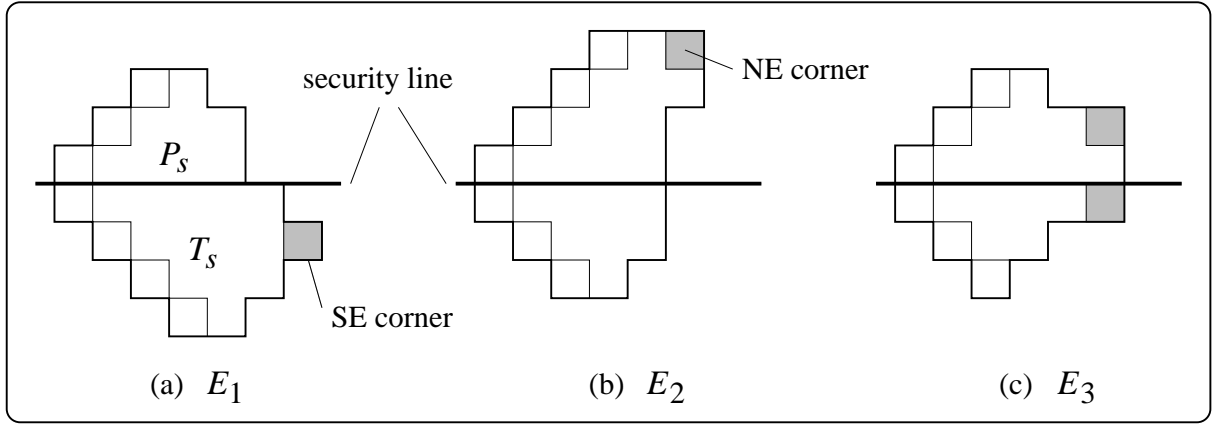}
\caption{\small The three classes of even doubly shifted stacks.}
\label{e1e2e3}
\end{center}
\end{figure}

We now concentrate on the computation of the generating series $E(x,z,w,q)$ and $A(x,z,w,q)$ of doubly shifted stacks. The concepts of flotation line, security line, North-East and South-East corners are defined similarly to the case of $D_S$-polyominoes. In the even case, we then introduce three subclasses of doubly shifted stacks as follows: $E_1$ if the South-East corner is below the security line; $E_2$ when the North-East corner is above the security line; and $E_3$ when both conditions are satisfied (see Figure~\ref{e1e2e3}). We obviously have
\begin{equation}
E(x,z,w,q) = E_1(x,z,w,q) + E_2(x,z,w,q) - E_3(x,z,w,q)\,.
\end{equation}
 
Observe that under the security line of an $E_1$-polyomino $Q$ lies a (single) shifted stack polyomino $T_S$ whereas on top of the security line is found a non-empty partition $P_S$ in distinct parts (rows). If $T_S$ has top width equal to $m$, then $P_S$ has biggest row at most $m$.

Let $T_{S,m}(x,y,q) = [u^m]T_S(u,x,y,q)$ be the generating series for shifted stacks having top width equal to $m$ (see Subsection~\ref{sec:stack}). We then have

\begin{prop}The generating series for even doubly shifted stacks $E_1$, $E_2$ and $E_3$ are given by
\begin{eqnarray}
E_1(x,z,w,q) & = & \sum_{m \geq 1}T_{S,m}(x,z,q)((-wq)_m) - 1)\,, \label{eqn:E1}\\
& & \nonumber \\
E_2(x,z,w,q) & = & E_1(x,w,z,q)\,, \label{eqn:E2} \\
& & \nonumber \\
E_3(x,z,w,q) & = & \sum_{m \geq 1}x^mzwq^{2m}(-zq)_{m-1}(-wq)_{m-1}\,. \label{eqn:E3}
\end{eqnarray}
\label{propA}
\end{prop}

\begin{proof}Formula~(\ref{eqn:E1}) follows from the previous discussion. To see~(\ref{eqn:E2}), simply observe that an horizontal reflection will send any $E_2$-polyomino into an $E_1$-polyomino and vice-versa. Finally, the closed form~(\ref{eqn:E3}) for $E_3(x,z,w,q)$ is obvious by construction.
\end{proof}

Note that a more closed form can be given for $E_1$, using the following functional equation. Let the variable $u$ mark the top row width of any convex polyomino and let $E_1(u) = E_1(u,x,z,w,q)$ and $T_S(u) = T_S(u,x,z,q)$ be the corresponding generating series.

\begin{prop}The generating series $E_1(u)$ of even doubly shifted stacks of the first class satisfies the functional equation
\begin{equation}
\label{E1funceqn}
E_1(u) = \frac{uqw}{1-uq}\left(T_S(1) - T_S(uq)\right) + \frac{w}{1-uq}\left(uqE_1(1) - E_1(uq)\right)\,.
\end{equation}
\label{propB}
\end{prop}

\begin{proof}The first term on the right hand side of~(\ref{E1funceqn}) has the effect of constructing the first row of the $d_2$-diagonal, directly above the base $T_S$, while the second term iteratively constructs the other rows, following the diagonal.
\end{proof}

\begin{figure}[!ht]
\begin{center}
\includegraphics{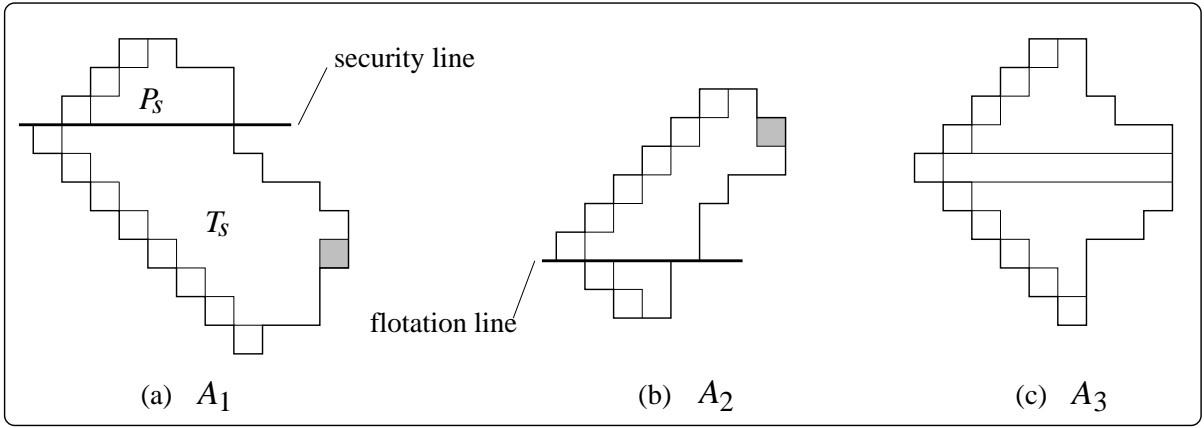}
\caption{\small The three classes acute doubly shifted stacks.}
\label{a1a2a3}
\end{center}
\end{figure}

In the odd case, we similarly introduce the following three subclasses of acute doubly shifted stacks: $A_1$ if the South-East corner is below the security line; $A_2$ if the North-East corner is above the flotation line; and $A_3$ if both conditions hold (see Figure~\ref{a1a2a3}). We obviously have
\begin{equation}
A(x,z,w,q) = A_1(x,z,w,q) + A_2(x,z,w,q) - A_3(x,z,w,q)
\end{equation}

The following results are similar to Propositions~\ref{propA} and~\ref{propB} above. Proofs are left to the reader.

\begin{prop}The generating series for acute doubly shifted stacks $A_1$, $A_2$ and $A_3$ are given by
\begin{eqnarray}
A_1(x,z,w,q) & = & \sum_{m \geq 1}T_{S,m}(x,z,q)(-wq)_{m-1}\,,\\
& & \nonumber \\
A_2(x,z,w,q) & = & A_1(x,w,z,q)\,,\\
& & \nonumber \\
A_3(x,z,w,q) & = & \sum_{m \geq 1}x^mzwq^m(-zq)_{m-1}(-wq)_{m-1}\,. \label{eqn:a3}
\end{eqnarray}\hfill $\square$
\end{prop}

\begin{prop}The generating series $A_1(u) = A_1(u,x,z,w,q)$ of acute doubly shifted stacks of the first class satisfies the functional equation
\begin{equation}
A_1(u) = wT_S(u,x,z,q) + \frac{w}{1-uq}\left(uqA_1(1) - A_1(uq)\right).
\end{equation}\hfill $\square$
\end{prop}

Finally, putting together formulas~(\ref{eqn:d1d2})--(\ref{eqn:a3}) yields $F_{d_1d_2}(t,q)$. Here are the first few terms:
\begin{equation}
F_{d_1d_2}(t,q) = t^{2}q+t^{4}q^{4}+t^{6}q^{5}+2\,t^{6}q^{7}+t^{6}q^{9}+2\,t^{8}q^{8}
+2\,t^{8}q^{10} + \ldots
\end{equation}

%
%
\section{Congruence-type convex polyominoes}

We are now in a position to enumerate congruence-type convex polyominoes, i.e. convex polyominoes up to rotation and reflection using formula~(\ref{BurnsideL}) of Section~\ref{sec:mobius} with $G = \mathfrak{D}_4$ and $\mathcal{P} = C$, the class of all convex polyominoes:
\begin{equation}
|C/\mathfrak{D}_4|_w = \frac{1}{8}\sum_{g \in \mathfrak{D}_4} |\Fix(g)|_w\,,
\end{equation}
where $|\Fix(g)|_w$ is the half-perimeter and area generating series of the convex $g$-symmetric polyominoes. Therefore,

\begin{prop}The half-perimeter and area generating series $(C/\mathfrak{D}_4)(t,q)$ of congruence-type polyominoes is given by
\begin{equation}
(C/\mathfrak{D}_4)(t,q) = |C/\mathfrak{D}_4|_w = \frac{1}{8} \left(C(t,t,q) + 2F_r(t,q) + F_{r^2}(t,t,q) + 2F_d(t,q)
+ 2F_v(t,q) \right)\,.
\end{equation}\hfill $\square$
\end{prop}

The first few terms of $(C/\mathfrak{D}_4)(t,q)$ are given in Appendix B.

%
%
\section{Asymmetric convex polyominoes}

We can also enumerate asymmetric convex polyominoes, i.e. polyominoes having no symmetry at all, using formula~(\ref{eqn:assym}) of Section~\ref{sec:mobius} and the value of the M\"obius function for the subgroups of $\mathfrak{D}_4$ given in Figure~\ref{TreillisIntro}.

\begin{prop}The half-perimeter and area generating series $\overline{C}(t,q) = F_{=0}$ of asymmetric convex polyominoes is given by
\begin{eqnarray}
\overline{C}(t,q) & = & F_{\geq 0} - F_{\geq d_1} - F_{\geq d_2} - F_{\geq r^2} - F_{\geq h} - F_{\geq v} +
2F_{\geq \langle d_1, d_2\rangle} + 2F_{\geq \langle h,v\rangle}\nonumber \\*
& & \nonumber \\*
& = & C(t,t,q) - 2F_d(t, q)  -  F_{r^2}(t,t,q) - 2F_v(t,q) + 2F_{d_1d_2}(t,q) + 2F_{hv}(t,q)\,.\label{eqn:asym} 
\end{eqnarray}
\end{prop}

The first few terms of this series are given in Appendix B.

Note that the same method would allow us to enumerate the convex polyominoes having exactly the symmetries of any given subgroup of $\mathfrak{D}_4$.

As an application we show the asymptotic result that for large area, almost all convex polyominoes are asymmetric. In other words, the probability for a convex polyomino to have at least one symmetry goes to zero as the area goes to infinity.

\begin{prop}\textnormal{(\textbf{Bender} \cite{EAB})}\ \ 
Let $c_n$ be the number of convex polyominoes with area $n$. Then
\begin{equation}
c_n \sim k\,\mu^n\,, \label{eqn:asympt}
\end{equation}
with
\begin{displaymath}
k = 2.67574\ldots \qquad \mu = 2.30913859330\ldots 
\end{displaymath}\hfill $\square$
\end{prop}

\begin{lemma}
Let $H$ be any non-trivial subgroup of $\mathfrak{D}_4$ and denote by $F_{H,n}$ the number of $H$-symmetric polyominoes with area $n$. Then,
\begin{equation}
\lim_{n \rightarrow \infty}\frac{F_{H,n}}{c_n} = 0\,.
\end{equation}
\end{lemma}

\begin{proof}A same basic argument works for $r^2$-, $h$-, $v$-, $d_1$- and $d_2$-symmetric convex polyominoes, using the fact that they are all constructed from two congruent subpolyominoes. A supplementary column or row has to be added in some cases, according to whether the height or the width of the initial polyomino is odd or even. These subpolyominoes are convex in every case, and they have at most half the area of the initial object.
\begin{itemize}
\item \textbf{$r^2$-symmetric convex polyominoes\,:}\ An $r^2$-symmetric convex polyomino with even width and area $n$ is constructed from two congruent convex subpolyominoes with exactly half the area, which can be glued together in at most $n/2$ ways (the maximal height of the columns that get glued together). So $F_{r^2,n}^{even} \leq \frac{1}{2}nc_{n/2}$. Hence
\begin{displaymath}
\lim_{n \rightarrow \infty}\frac{F_{r^2,n}^{even}}{c_n} \quad \leq \quad  \lim_{n \rightarrow \infty}\frac{\frac{1}{2}nk\mu^{n/2}}{k\mu^n} \quad = \quad 0\,.
\end{displaymath}

Next consider an $r^2$-symmetric convex polyomino with odd width and area $n$. This polyomino is constructed from a central column ($n$ choices of height) and two congruent convex subpolyominoes of area at most $\lfloor n/2 \rfloor$. Then there are at most $n$ possible positions where to glue the subpolyominoes to the central column (they are glued symmetrically). Thus $F_{r^2,n}^{odd} \leq n^2(1+c_1+c_2+\ldots + c_{\lfloor n/2 \rfloor}) < n^3c_{\lfloor n/2 \rfloor}$ and the result follows as above. Hence the result holds for the subgroup $\langle r^2 \rangle$ of $\mathfrak{D}_4$;

\item \textbf{$h$- and $v$-symmetric convex polyominoes\,:}\ A $v$-symmetric convex polyomino of even width is constructed much in the same way as an $r^2$-symmetric convex polyomino of even width, with the difference that there is only one glueing position for the two subpolyominoes. So $F_{v,n}^{even} \leq c_{n/2}$ and the limit will be zero. Similarly, the construction for $v$-symmetric convex polyominoes of odd width resembles that of $r^2$-symmetric convex polyominoes of odd width, and we get $F_{v,n}^{odd} \leq n^2(1+c_1+c_2+\ldots + c_{\lfloor n/2 \rfloor}) < n^3c_{\lfloor n/2 \rfloor}$. So the result holds for the subgroups $\langle h \rangle$ and $\langle v \rangle$;

\item \textbf{$d_1$- and $d_2$-symmetric convex polyominoes\,:}\ Let $P$ be a $d$-symmetric convex polyomino and $Q$ its fundamental region. Suppose that $P$ has $b$ cells on the diagonal symmetry axis. Then the minimum area $P$ can have is $b + 2(b-1)$. This gives a minimum area of $b + (b-1)$ for $Q$. Hence
\begin{displaymath}
\frac{\textrm{Area of } Q_{min}}{\textrm{Area of } P_{min}} \quad = \quad \frac{2b-1}{3b-2}\,. 
\end{displaymath}

Then if we add to $Q$ a cell not on the diagonal symmetry axis, two cells get added to $P$, and thus we conclude that the ratio can only decrease as we make $P$ into a larger $d$-symmetric convex polyomino with the same number of cells on the diagonal axis. For $b \geq 2$, the ratio will be smaller than or equal to $3/4$. As a loose approximation, we can take $Q$ to be any convex polyomino. This gives $F_{d,n} \leq c_{\lceil 3n/4 \rceil} + c_{\lceil n/2 \rceil}$. The $c_{\lceil n/2 \rceil}$ term corresponds to those $d$-symmetric convex polyominoes having only one cell on the diagonal axis (in which case the fundamental region is a convex polyomino of area $(n+1)/2$). Hence the result will also hold for the subgroups $\langle d_1 \rangle$ and $\langle d_2 \rangle$. 
\end{itemize}
 
Then since $r$-symmetric polyominoes are a subclass of $r^2$-symmetric polyominoes, $hv$-symmetric polyominoes a subclass of $h$- or $v$-symmetric polyominoes and $d_1d_2$-symmetric polyominoes a subclass of $d_1$- or $d_2$-symmetric polyominoes, the results extends to all of the symmetry classes we considered.
\end{proof}

\begin{prop}
If we denote by $\overline{c}_{n}$ the number of asymmetric convex polyominoes of area $n$, then
\begin{equation}
\overline{c}_{n} \sim c_n\,.
\end{equation}
\end{prop}

\begin{proof}
 We get the result from equation~(\ref{eqn:asym}) and the previous lemma.
\end{proof}


\pagebreak

\noindent{\LARGE Appendix A}
\vspace{5mm}

The following two tables present the numbers of convex polyominoes according to their symmetry types and their perimeter or area. The columns indexed by subgroups of $\mathfrak{D}_4$ give the numbers of convex polyominoes of a given perimeter or area that are left fixed by the symmetries of the subgroup. The columns \textsl{rotation, congruence} and \textsl{asym} give respectively the number of rotation-type, congruence-type and asymmetric convex polyominoes of the given size.

\begin{table}[!ht]
\begin{center}
\begin{scriptsize}
\begin{tabular}{|c||c|c|c|c|c|c|c|c|c|c|}
\hline
per & \{1\} & $\langle r\rangle, \langle r^3\rangle$ & $\langle r^2\rangle$ & rotation & $\langle h\rangle, \langle v\rangle$ & $\langle d_1\rangle, \langle d_2\rangle$  & congruence & $\langle h,v\rangle$ & $\langle d_1,d_2\rangle$ & asym\\
\hline \hline
4  & 1     & 1 & 1   & 1     & 1   &1   & 1    & 1  &1  & 0    \\ \hline
6  & 2     & 0 & 2   & 1     & 2   &0   & 1    & 2  &0  & 0    \\ \hline
8  & 7     & 1 & 3   & 3     & 3   &3   & 3    & 3  &1  & 0    \\ \hline
10 & 28    & 0 & 8   & 9     & 6   &0   & 6    & 4  &0  & 16   \\ \hline
12 & 120   & 2 & 16  & 35    & 12  &14  & 24   & 6  &4  & 72   \\ \hline
14 & 528   & 0 & 40  & 142   & 24  &0   & 77   & 8  &0  & 456  \\ \hline
16 & 2344  & 4 & 84  & 609   & 48  &70  & 334  & 12 &8  & 2064 \\ \hline
18 & 10416 & 0 & 200 & 2654  & 96  &0   & 1351 & 16 &0  & 10056\\ \hline
20 & 46160 & 8 & 424 & 11650 & 192 &348 & 5960 & 24 &24 & 44736\\ \hline
\end{tabular}
\end{scriptsize}
\caption{\small Convex polyominoes enumerated by their symmetries and perimeter.}
\end{center}
\end{table}

\begin{table}[!ht]
\begin{center}
\begin{scriptsize}
\begin{tabular}{|c||c|c|c|c|c|c|c|c|c|c|}
\hline
area & \{1\} & $\langle r\rangle, \langle r^3\rangle$ & $\langle r^2\rangle$ & rotation & $\langle h\rangle, \langle v\rangle$ & $\langle d_1\rangle, \langle d_2\rangle$  & congruence & $\langle h,v\rangle$ & $\langle d_1,d_2\rangle$ & asym\\
\hline \hline
1  & 1    & 1 & 1   & 1    & 1  & 1  & 1    & 1 & 1 & 0    \\ \hline
2  & 2    & 0 & 2   & 1    & 2  & 0  & 1    & 2 & 0 & 0    \\ \hline
3  & 6    & 0 & 2   & 2    & 2  & 2  & 2    & 2 & 0 & 0    \\ \hline
4  & 19   & 1 & 7   & 7    & 5  & 1  & 5    & 3 & 1 & 8    \\ \hline
5  & 59   & 1 & 7   & 17   & 5  & 5  & 11   & 3 & 1 & 40   \\ \hline
6  & 176  & 0 & 24  & 50   & 12 & 4  & 29   & 4 & 0 & 128  \\ \hline
7  & 502  & 0 & 22  & 131  & 12 & 14 & 72   & 4 & 2 & 440  \\ \hline
8  & 1374 & 2 & 74  & 363  & 26 & 12 & 191  & 6 & 2 & 1240 \\ \hline
9  & 3630 & 2 & 62  & 924  & 26 & 38 & 478  & 6 & 2 & 3452 \\ \hline
10 & 9312 & 0 & 208 & 2380 & 52 & 32 & 1211 & 6 & 2 & 8952 \\ \hline
\end{tabular}
\end{scriptsize}
\caption{\small Convex polyominoes enumerated by their symmetries and area.}
\end{center}
\end{table}

\pagebreak

\noindent{\LARGE Appendix B}
\vspace{5mm}

We provide here the beginning terms of the series $(C/\mathfrak{C}_4)(t,q)$, $(C/\mathfrak{D}_4)(t,q)$ and $\overline{C}(t,q)$ of rotation-type, congruence-type and asymmetric polyominoes respectively. They are given up to area 10 ($q^{10}$).

\begin{eqnarray*}
(C/\mathfrak{C}_4)(t,q) & = & t^{2}q+t^{3}q^{2}+2\,t^{4}q^{3}+t^{4}q^{4}+6\,t^{5}q^{4}+2\,t^{5}q^{5}+t^{5}q^{6}+15\,t^{6}q^{5}+11\,t^{6}q^{6}+6\,t^{6}{q}^{7}\\
&& +38\,t^{7}q^{6}+2\,t^{6}q^{8}+38\,t^{7}q^{7}+t^{6}q^{9}+36\,t^{7}q^{8}+87\,t^{8}q^{7}+18\,t^{7}q^{9}+124\,t^{8}q^{8}+9\,t^{7}q^{10}\\
&& +139\,t^{8}q^{9}+201\,t^{9}q^{8}+115\,t^{8}q^{10}+334\,t^{9}q^{9}+470\,t^{9}q^{10}+432\,t^{10}q^{9}+861\,t^{10}{q}^{10}\\
&& +925\,t^{11}q^{10} + \ldots\\
&& \\
(C/\mathfrak{D}_4)(t,q) & = & {t}^{2}q+{t}^{3}{q}^{2}+2\,{t}^{4}{q}^{3}+{t}^{4}{q}^{4}+4\,{t}^{5}{q}^{4}+{t}^{5}{q}^{5}+{t}^{5}{q}^{6}+10\,{t}^{6}{q}^{5}+7\,{t}^{6}{q}^{6}+4\,{t}^{6}{q}^{7}\\
&& +21\,{t}^{7} {q}^{6}+2\,{t}^{6}{q}^{8}+19\,{t}^{7}{q}^{7}+{t}^{6}{q}^{9}+20\,{t}^{7}{q}^{8}+49\,{t}^{8}{q}^{7}+9\,{t}^{7}{q}^{9}+65\,{t}^{8}{q}^{8}+6\,{t}^{7}{q}^{10 }\\
&& +74\,{t}^{8}{q}^{9}+104\,{t}^{9}{q}^{8}+62\,{t}^{8}{q}^{10}+167\,{t}^{9}{q}^{9}+239\,{t}^{9}{q}^{10}+227\,{t}^{10}{q}^{9 }+436 \,{t}^{10}{q}^{10}\\
&& +468\,{t}^{11}{q}^{10} + \ldots\\
&& \\
\overline{C}(t,q) & = & 8\,{t}^{5}{q}^{4}+8\,{t}^{5}{q}^{5}+32\,{t}^{6}{q}^{5}+24\,{t}^{6}{q}^{6}+16\,{t}^{6}{q}^{7}+104\,{t}^{7}{q}^{6}+152\,{t}^{7}{q}^{7}+104\,{t}^{7}{q}^{8}\\
&& +272\,{t}^{8}{q}^{7}+72\,{t}^{7}{q}^{9}+448\,{t}^{8}{q}^ {8}+16\,{t}^{7}{q}^{10}+496\,{t}^{8}{q}^{9}+688\,{t}^{9}{q}^{8}+400\,{t}^{8}{q}^{10}\\
&& +1336\,{t}^{9}{q}^{9}+ 1744\,{t}^{9}{q}^{10}+1552\,{t}^{10}{q}^{9}+3344\,{t}^{10}{q}^{10}+3448\,{t}^{11}{q}^{10}+ \ldots
\end{eqnarray*}

\end{document}